\documentclass[a4paper,12pt,reqno]{amsart}

\usepackage{amssymb,amsmath,amsfonts,amsthm,mathtools,dsfont}
\usepackage{latexsym,mathrsfs,pb-diagram}
\usepackage[pdftex]{graphicx}
\usepackage{pgf,tikz,tikz-cd}

\usepackage[hmargin=3.0cm,vmargin=2.5cm]{geometry}
\usepackage[plainpages=false,colorlinks,pdfpagelabels]{hyperref}
\hypersetup{
	urlcolor=black,
	citecolor=green,
	linkcolor=blue
}

\numberwithin{equation}{section}

\theoremstyle{definition}
\newtheorem{Def}{Definition}[section]

\theoremstyle{remark}

\newtheorem{Rem}[Def]{Remark}

\theoremstyle{plain}
\newtheorem{Prop}[Def]{Proposition}
\newtheorem{Cor}[Def]{Corollary}
\newtheorem{Thm}[Def]{Theorem}
\newtheorem{Lem}[Def]{Lemma}

\newcommand{\Z}{\mathbb{Z}}
\newcommand{\N}{\mathbb{N}}
\newcommand{\R}{\mathbb{R}}

\newcommand{\C}{\mathbb{C}}

\DeclareMathOperator{\Span}{\mathrm{span}}

\newcommand{\dd}{\mathrm{d}}

\newcommand{\TT}{\mathbb{T}}

\newcommand{\D}{\mathscr{D}}

\newcommand{\cinfty}{\mathscr{C}^\infty}

\newcommand{\DD}{\mathrm{D}}

\DeclareMathOperator{\supp}{\mathrm{supp}}

\newcommand{\wf}{\mathfrak s}

\title{Propagation of regularity for a class of systems of real vector fields on torus}

\author{I. A.~Ferra}
\address{Universidade Federal do ABC, Brazil}
\email{\texttt{ferra.igor@ufabc.edu.br}}

\author{L. A.~Carvalho dos Santos}
\address{Universidade Federal de São Carlos, Brazil}
\email{\texttt{luiscarvalho@ufscar.br}}

\thanks{Supported by
	Fundação de Amparo à Pesquisa do Estado de São Paulo (FAPESP) under grant 2023/11769-5 and by
	Conselho Nacional de Desenvolvimento Cient{\'i}fico e Tecnol{\'o}gico
	(CNPq) under  grant 404175/2023-6.
}

\begin{document}

\maketitle

\begin{abstract}
	We characterize the global hypoellipticity, almost hypoellipticity and solvability for a class of systems of real vector fields on the $(n+1)$-dimensional torus as well as the same properties about the sum of squares associated to the system. The key result is a theorem about propagation of regularity for solutions of the system.
\end{abstract}

\section{Introduction}
A smooth real $1$-form $a(t)=\sum_{j=1}^n a_j(t)\dd t_j$ in the $n$-dimensional torus $\TT^n$, where $(t,x)=\left(t_1,\ldots,t_n,x\right)$ denotes standard angular coordinates on $\TT^{n+1}$, defines a system $\mathcal X$ of differential operators given by the real vector fields
\begin{align}
	\label{eq:intro-system}
X_j = \partial_{t_j} + a_j(t)\partial_x,\quad j \in \{1,\ldots,n\}.
\end{align}
When $a$ is closed, this kind of system is associated to a differential complex that comes from a locally integrable structure (see \cite{bcm93, bch_iis}).
The closedness of $a$ is equivalent to the commutativity of the operators $X_1,\ldots,X_n$ and, in this case, standard arguments can be used in order to characterize the global hypoellipticity and solvablity of $\mathcal X$ in terms of diophantine conditions (see \cite{bcm93,bkmz16,bp99sys,hz19}). The topic of regularity of the solutions as well as the solvability for systems of (not necessarily real) vector fields or even pseudodifferential operators that commute, on torus, has been widely studied in recent years \cite{akm19, bk07, bkmz16, bmz12, af21}.

One of the main goals of this work is to obtain characterizations of global hypoellipticity or solvability for the non-commutative case. Removing the commutativity hypothesis may seem artificial from the point of view of involutive structure theory, but the same does not occur if we investigate problems about the sum of squares $P=\sum_{j=1}^n X_j^2$. Our goal is then justified by the following
\begin{Thm}
	\label{Thm:P-GH-iff-L-GH}
	Let $u \in \D'\left(\TT^{n+1}\right)$. Then $Pu \in \cinfty\left(\TT^{n+1}\right)$ if and only if $X_ju \in \cinfty\left(\TT^{n+1}\right)$ for every $j\in\{1,\ldots,n\}$. In particular, the system $\mathcal X$ is globally hypoelliptic in $\TT^{n+1}$ if and only if the operator $P$ is globally hypoelliptic in $\TT^{n+1}$.
\end{Thm}
%



 On torus, this result was known only for operators with constant coefficients (see Proposition 3.11 in~\cite{afr22}). The global hypoellipticity for $P$ when $n=2$ was studied in~\cite{hps06,hps12} where it was proved that $P$ is globally hypoelliptic if and only if there is a point of finite type for the vector fields $X_1$ and $X_2$ (i.e. the H\"ormander's condition is satisfied in a point) or some diophantine condition holds true. For an arbitrary $n$, it is not difficult to see that the system $\mathcal X$ is commutative if and only if the H\"ormander's condition fails in every point for $X_1,\ldots,X_n$. Unfortunately, the technique in~\cite{hps06,hps12} does not seem to apply to higher-dimensional tori.

Regarding hypoellipticity of $P$ in the general case (arbitrary $n$), if $\mathcal X$ is not commutative then $P$ is locally hypoelliptic in some open subset $U\times \TT\subset \TT^n\times \TT$ (by the H\"ormander's condition) and, therefore, the study of propagation of regularity turns out to be a natural approach to study the global hypoellipticity of $P$. Indeed, we prove in Section~\ref{sec:PR} that in the non-commutative case, the operator $P$ is globally	hypoelliptic by showing the following strong result about propagation of regularity (recall that by Theorem \ref{Thm:P-GH-iff-L-GH} the conditions $Pu\in\cinfty\left(\TT^{n+1}\right)$ and $X_j u \in \cinfty\left(\TT^{n+1}\right)$ for every $j\in\{1,\ldots,n\}$ are equivalent):  
\begin{Thm}
	\label{Thm:propagation-one-point}
	Suppose that $X_ju \in \cinfty\left(\TT^{n+1}\right)$ for every $j\in\{1,\ldots,n\}$ and that there exists $s=\left(s_1,\ldots,s_n\right)\in\TT^n$ such that given $k\in\Z_+$ we find $C_k>0$ such that
	\[
	\left|\hat u(s,\xi)\right| \le C_k\left(1+|\xi|\right)^{-k}, \quad \forall \xi \in \Z.
	\]
	Then $u \in \cinfty\left(\TT^{n+1}\right)$.
\end{Thm}

 Results of this nature were obtained only for involutive systems (see the proof of Theorem 2.4 in~\cite{bcm93}). For sums of squares, previous works (see \cite{afr22,hp02}) provided results about propagation of regularity only when the sum of squares has the Laplacian of the $t$-variable in its expression. Therefore, Theorem \ref{Thm:propagation-one-point} is a novelty for both systems and sums of squares of real vector fields. We also stress that the previous results about propagation of regularity -- as well more results about the global hypoellipticity for sums of squares (see \cite{albanese11, afr24, bfp17, hp00}) --  use energy $L^2$-estimates in its proof, which is greatly facilitated in the presence of the $t$-Laplacian term. The second great novelty of Theorem \ref{Thm:propagation-one-point} is the use of a new and simpler thecnique that makes slightly more accurate estimates without using $L^2$-norms. It is worth to point out that in~\cite{hps06, hps12}, the authors first exhaust the study about diophantine conditions and then they look at the finite type condition. Here, Theorem \ref{Thm:propagation-one-point} allows us to do the opposite: we first obtain the result when there is a point of finite type and, when all the points are of infinite type, $a$ is closed and we can just use the vast literature on the subject to obtain characterization of global hypoellipticity in terms of a diophantine condition. This approach, in addition to generalizing the results in~\cite{hps06,hps12}, considerably simplifies the proof given there.

Furthermore, since global hypoellipticity implies almost global hypoellipticity and global solvability (see~\cite{afr22} and also see Section~\ref{sec:global-solv} for the precise definition of these concepts), when the system is not commutative these three properties hold true. In particular, our results also allow us to give a complete characterization of these properties for $\mathcal X$ and $P$ (with a simple proof). Also, from our results we conclude that the most interesting case to study is the well known case that comes from an involutive structure. This is the case where a diophantine condition plays an important role and to state our generalization of the results in~\cite{hps06,hps12} we denote
\begin{align}
	\label{eq:averages}
a_{j0}(t) = \frac{1}{2\pi} \int_0^{2\pi} a_{j}\left(t_1,\ldots,t_{j-1},\lambda,t_{j+1},\ldots,t_n\right)\dd \lambda,\quad t \in \TT^n, j \in\{1,\ldots,n\}.
\end{align}
\begin{Thm}
	\label{Thm:equivalence-gh}
The operator $P$ is globally hypoelliptic in $\TT^{n+1}$ if and only if there is a point $t \in \TT^n$ of finite type for $X_1,\ldots,X_n$ or the range of $a_{j0}$ contains a real number $\alpha_j$ for each $j\in\{1,\ldots,n\}$ such that the collection $\{\alpha_1,\ldots,\alpha_n\} \subset \R$ is not simultaneously approximable.
\end{Thm}
The proof of this result, as well as the definition of simultaneously approximability, are offered in Section~\ref{sec:GH}. We also follow the ideas in~\cite{afr24} and replace these diophantine conditions by the global hypoellipticity of a system of real vector fields with constant coefficients given by
\[
\mathcal X_0 = \left\{\partial_{t_{j}} + a_{j0}\left(t\right)\partial_x : j \in \{1,\ldots,n\}, t \in \TT^n\right\}.
\]
For the global hypoelipticity we prove 
\begin{Thm}
	\label{Thm:equivalence-GH-sums-systems}
	The operator $P$ is globally hypoelliptic in $\TT^{n+1}$ if and only if $\mathcal X_0$ is globally hypoelliptic in $\TT^{n+1}$ or there exists a point $t \in \TT^n$ of finite type for $X_1,\ldots,X_n$.
\end{Thm}
For the almost global hypoellipticity and global solvability we have
\begin{Thm}
	\label{Thm:most-general-equiv}
	The following properties are equivalent:
	\begin{enumerate}
		\item $P$ is globally solvable;
		\item $P$ is almost globally hypoelliptic;
		\item $\mathcal X$ is globally solvable;
		\item $\mathcal X$ is almost globally hypoelliptic;
		\item $\mathcal X_0$ is globally solvable or there is a point of finite type for $X_1,\ldots,X_n$;
		\item $\mathcal X_0$ is almost globally hypoelliptic or there is a point of finite type for $X_1,\ldots,X_n$.
	\end{enumerate}
\end{Thm}
This work is organized as follows: we first prove the central result about propagation of regularity, Theorem \ref{Thm:propagation-one-point}, in Section~\ref{sec:PR}. In order to prove Theorem~\ref{Thm:P-GH-iff-L-GH}, which we do in Section~\ref{sec:system-sums}, we need some microlocal results developed in Appendix~\ref{sec:me}. Although Corollary \ref{local_ellipticity_su} and Corollary \ref{Cor:ellipticity-wf} are classic results of microlocal regularity, we choose to keep them in the text since the version we need is about the set $\mathfrak{s}(u)$ for some $u \in \mathscr{D}'(\TT^{n+1})$ (see Appendix~\ref{sec:me} for the definition of $\wf(u)$) rather than the usual wave front set of $u$. 
In Section \ref{sec:GH} we prove Theorem~\ref{Thm:equivalence-gh} and Theorem~\ref{Thm:equivalence-GH-sums-systems} and we close the work with Section~\ref{sec:global-solv}, with the proof of Theorem~\ref{Thm:most-general-equiv}.

\section{Propagation of regularity}
\label{sec:PR}

In this section we prove our main theorem about propagation of regularity of the system $\mathcal X$ defined by \eqref{eq:intro-system}. For each $j \in \{1,\ldots,n\}$ we define
\[
A_j:\TT^n\to \R, \,\, A_j(t)=\int_0^{t_j} a_{j}\left(t_1,\ldots,t_{j-1},\lambda,t_{j+1},\ldots,t_n\right)\dd \lambda - a_{j0}(t)t_j \in \cinfty\left(\TT^n\right),
\]
where $a_{j0}$ is given by \eqref{eq:averages} (this function $a_{j0}$ belongs to $\cinfty\left(\TT^n\right)$ and does not depend on $t_j$). We can use $A_j$ to set a diffeomorphism
\[
\psi_j : \TT^n\times \TT \to \TT^n\times \TT, \,\, \psi_j\left(t,x\right) = \left(t, x+A_j(t)\right),
\]
which in turn induces an automorphism of $\D'\left(\TT^n\times \TT\right)$ that restricts itself to an automorphism of $\cinfty\left(\TT^n\times\TT\right)$ given by
\[
S_j : \D'\left(\TT^n\times \TT\right)\to \D'\left(\TT^n\times \TT\right), \,\, S_j u \doteq \sum_{\xi \in \Z} e^{iA_j(t)\xi} \hat u(t,\xi) e^{ix\xi}.
\]
%
So $S_ju$ is characterized by the relation
\[
\widehat{S_ju}(t,\xi) = e^{iA_j(t)\xi} \hat u(t,\xi) \in \D'\left(\TT^n\right), \quad \forall \xi \in \Z.
\]
Since each $A_{j}$ is a real function, if $\hat u\left(\cdot,\xi\right) \in \cinfty\left(\TT^n\right)$ then
\begin{align}
	\label{eq:partial-fourier-auto-Sj-modulos}
	\left|\widehat{S_ju}(t,\xi)\right| = \left|\hat u(t,\xi)\right|, \quad \forall t\in\TT^n, j\in\{1,\ldots,n\}.
\end{align}
Note also that if $Y_j\doteq \partial_{t_j}+ a_{j0}(t)\partial_{x}$ and $u\in \D'\left(\TT^{n+1}\right)$ then
\begin{align*}
	Y_j\left(S_j u\right) &= \sum_{\xi \in \Z} \left(\partial_{t_j}+a_{j0}(t)\partial_{x}\right)e^{iA_j(t)\xi} \hat u(t,\xi) e^{ix\xi}\\
	&= \sum_{\xi \in \Z} e^{iA_j(t)\xi} \left(\partial_{t_j}+i\partial_{t_j}\left(A_j(t)\xi\right)+ia_{j0}(t)\xi\right)\hat u(t,\xi) e^{ix\xi}\\
	&= \sum_{\xi \in \Z} e^{iA_j(t)\xi} \left(\partial_{t_j}+i \left(a_{j}(t)-a_{j0}(t)\right)\xi+ia_{j0}(t)\xi\right)\hat u(t,\xi) e^{ix\xi}\\
	&= \sum_{\xi \in \Z} e^{iA_j(t)\xi} \left(\partial_{t_j}+ia_{j}(t)\xi\right)\hat u(t,\xi) e^{ix\xi}\\
	&= S_j\left(\sum_{\xi \in \Z} \left(\partial_{t_j}+ia_{j}(t)\xi\right)\hat u(t,\xi) e^{ix\xi}\right)\\
	&= S_j\left(\partial_{t_j}+a_{j}(t)\partial_{x}\right)u\\
	&= S_j X_j u.
\end{align*}
Taking the Fourier transform with respect to $x$,
\begin{align}
	\label{eq:automorphism_Sj-fourier-formula}
	\left(\partial_{t_j}+ia_{j0}(t)\xi\right)\widehat{S_j u}(t,\xi) = \widehat{S_jX_j u}(t,\xi), \quad \forall \xi \in \Z.
\end{align}
\begin{Rem}
	\label{Rem:infinite-type}
	Let us suppose that $\partial_{t_j}a_k - \partial_{t_k}a_j = 0$, that is the case where there is no point of finite type for $X_1,\ldots,X_n$. It follows that each $a_{j0}$ is constant: if $\ell \ne j$ then
	\[
	2\pi \partial_{t_\ell} a_{j0}\left(t\right) = \int_0^{2\pi} \partial_{t_j} a_{\ell}\left(t_1,\ldots,t_{j-1},\lambda,t_{j+1},\ldots,t_n\right)\dd\lambda = 0.
	\]
	So we can write $a_{j0}(t) = \alpha_j \in \R$ for every $t \in \TT^n$. Also in this case (see~\cite{bp99sys}) there is $A \in \cinfty\left(\TT^n\right)$ such that
	\[
	\partial_{t_j}A(t) = a_j(t)-\alpha_j, \quad \forall j \in\{1,\ldots,n\}, t \in \TT^n.
	\] 
	Indeed, take
	\[
	A(t) = \sum_{j=1}^n \left(\int_{0}^{t_j} a_{j}\left(0,\ldots,\lambda,\ldots,t_n\right)\dd \lambda - \alpha_j t_j\right),
	\]
	where $\lambda$ is in the $j$-th position, which belongs to $\cinfty\left(\TT^n\right)$ since $a_{j0}$ is constant for every $j\in\{1,\ldots,n\}$. Note that the $j$-th term of the sum above depends only on the variables $t_{j},\ldots,t_n$, so
	\[
	\partial_{t_1} A(t) = a_1(t)-\alpha_1,
	\]
	\begin{align*}
		\partial_{t_2} A(t) &= \partial_{t_2}\sum_{j=1}^2\left(\int_{0}^{t_j} a_{j}\left(0,\ldots,\lambda,\ldots,t_n\right)\dd \lambda - \alpha_j t_j\right)\\
		&= \int_0^{t_1} \partial_2 a_1\left(\lambda,t_2,\ldots,t_n\right)\dd \lambda + a_2\left(0,t_2,\ldots,t_n\right) - \alpha_2\\
		&= \int_0^{t_1} \partial_1 a_2\left(\lambda,t_2,\ldots,t_n\right)\dd \lambda + a_2\left(0,t_2,\ldots,t_n\right) - \alpha_2\\
		&= a_2\left(t\right) - \alpha_2,
	\end{align*}
	\begin{align*}
		\partial_{t_3} A(t) &= \partial_{t_3}\sum_{j=1}^3\left(\int_{0}^{t_j} a_{j}\left(0,\ldots,\lambda,\ldots,t_n\right)\dd \lambda - \alpha_j t_j\right)\\
		&= \int_0^{t_1} \partial_3 a_1\left(\lambda,t_2,\ldots,t_n\right)\dd \lambda + \int_0^{t_2} \partial_3 a_2\left(0,\lambda,\ldots,t_n\right)\dd \lambda+ a_3\left(0,0,t_3,\ldots,t_n\right) - \alpha_3\\
		&= \int_0^{t_1} \partial_1 a_3\left(\lambda,t_2,\ldots,t_n\right)\dd \lambda + \int_0^{t_2} \partial_2 a_3\left(0,\lambda,\ldots,t_n\right)\dd \lambda + a_3\left(0,0,t_3,\ldots,t_n\right) - \alpha_3\\
		&= a_3\left(t\right) - \alpha_3
	\end{align*}
	and so on. If we set
	\[
	S:\D'\left(\TT^{n+1}\right) \to \D'\left(\TT^{n+1}\right), \,\, Su=\sum_{\xi \in \Z} \hat u(t,\xi) e^{iA(t)\xi}e^{ix\xi}
	\]
	then as we have done for \eqref{eq:automorphism_Sj-fourier-formula} we can demonstrate that
	\begin{align}
		\label{eq:simultaneous-reduction}
		Y_j S = SX_j, \quad \forall j \in \{1,\ldots,n\}.
	\end{align}
	In other words, when the H\"ormander's condition fails everywhere, there is a conjugation between the original system and a system of operators with constant coefficients.
\end{Rem}

In order to prove Theorem \ref{Thm:propagation-one-point}, we need some auxiliary results. We start with the following
\begin{Lem}
	\label{Lem:estimative-zeroth-derivative}
	Suppose that $u\in \D'\left(\TT^{n+1}\right)$ is such that $X_j u \in \cinfty\left(\TT^{n+1}\right)$ for every $j\in\{1,\ldots,n\}$. Then $u \in \cinfty\left(\TT^{n+1}\right)$ provided given $k\in\Z_+$ there exists $C_k>0$ such that
	\[
	\left|\hat u(t,\xi)\right| \le C_k\left(1+|\xi|\right)^{-k}, \quad \forall t \in \TT^n, \xi \in \Z.
	\]
\end{Lem}
\begin{proof}
	First notice that the hypothesis ensures that $Pu \in \cinfty\left(\TT^{n+1}\right)$, so by ellipticity we conclude that $\hat u(\cdot,\xi) \in \cinfty\left(\TT^n\right)$ for every $\xi \in \Z$. To prove that $u\in \cinfty\left(\TT^{n+1}\right)$ it is enough to show that given $k\in\Z_+$ and $\alpha\in\Z^n_+$, there exists $C_{\alpha,k}>0$ such that
	\begin{align}
		\label{eq:estimative-smooth-decay-alpha-N}
		\left|\partial_t^\alpha \hat u(t,\xi)\right| \le C_{\alpha,k}\left(1+|\xi|\right)^{-k},\quad \forall t\in \TT^n, \xi \in \Z.
	\end{align}
	This will be done by induction on $|\alpha|$. For $|\alpha|=0$ we can take $C_{0,k} = C_k$. Suppose that \eqref{eq:estimative-smooth-decay-alpha-N} is valid for $|\alpha|\le \ell$ and take $\beta \in \Z_+^{n}$ such that $|\beta|=\ell+1$. If we write $\beta=\left(\beta_1,\ldots,\beta_n\right)$ then there exists $j\in\{1,\ldots,n\}$ such that $\beta_j\ne 0$ and we set $\alpha=\beta-e_j$, where $e_j$ denotes the $j$-th vector of the canonical basis of $\R^n$. By the inductive hypothesis given $k\in\Z_+$, there is $B_{\ell,k}>0$ such that
	\[
	\left|\partial_t^\gamma \hat u(t,\xi)\right| \le B_{\ell,k}\left(1+|\xi|\right)^{-k},\quad \forall t\in \TT^n, \xi \in \Z, |\gamma|\le \ell.
	\]
	Then 
	\begin{align*}
		\left|\partial_t^\gamma \widehat{\partial_{x}u}(t,\xi)\right| &\le |\xi| \left|\partial_t^\gamma \hat u(t,\xi)\right| \\
		&\le B_{\ell,k+1}\left(1+|\xi|\right)^{-k},\quad \forall t\in \TT^n, \xi \in \Z, |\gamma|\le \ell.
	\end{align*}
	If we apply $\partial_t^\alpha$ in the equality 
	\[
	\partial_{t_j} \hat u(t,\xi) = \widehat{X_ju}(t,\xi) - a_j(t)\widehat{\partial_xu}(t,\xi)
	\]
	and use that there is $C_{\alpha,j,k}>0$ such that
	\[
	\left|\partial_t^\alpha \widehat{X_ju}(t,\xi)\right|\le C_{\alpha,j,k}\left(1+|\xi|\right)^{-k},\quad \forall (t,\xi)\in \TT^n\times \Z, j\in\{1,\ldots,n\}
	\]
	we obtain that
	\begin{align*}
		\left|\partial_t^\beta \hat u(t,\xi)\right| & = \left|\partial_t^\alpha\widehat{\partial_{t_j} u}(t,\xi)\right| \\
		&\le \left|\partial_t^\alpha \widehat{X_ju}(t,\xi)\right| + \left|\sum_{\gamma\le\alpha}\binom{\alpha}{\gamma}\left[\partial_t^\gamma \widehat{\partial_x u}(t,\xi)\right] \partial_t^{\alpha-\gamma}a_j(t)\right|\\
		&\le C_{\alpha,j,k}\left(1+|\xi|\right)^{-k} +  \sum_{\gamma\le\alpha}\binom{\alpha}{\gamma} \left\|\partial_t^{\alpha-\gamma}a_j\right\|_{L^\infty\left(\TT^n\right)}  B_{\ell,k+1}\left(1+|\xi|\right)^{-k}\\
		&= C_{\beta,k} \left(1+|\xi|\right)^{-k}, \quad \forall t \in \TT^n, \xi \in \Z,
	\end{align*}
	where
	\[
	C_{\beta,k} = \max_{1\le j \le n} C_{\alpha,j,k} + B_{\ell,k+1}\sum_{\gamma\le\alpha}\binom{\alpha}{\gamma} \left\|\partial_t^{\alpha-\gamma}a_j\right\|_{L^\infty\left(\TT^n\right)} .
	\]
\end{proof}
%
%
%
\begin{Lem}
	\label{Lem:sobolev-estimate}
	Suppose that $u\in \D'\left(\TT^{n+1}\right)$ and $Pu \in \cinfty\left(\TT^{n+1}\right)$. Then there are $D>0$ and $M\in\Z_+$ such that
	\[
	\left|\hat u(t,\xi)\right| \le D\left(1+|\xi|\right)^M, \quad \forall t \in \TT^n, \xi \in \Z.
	\]
\end{Lem}
\begin{proof}
	First we take $p \in \Z_+$ such that $u \in H^{-p}\left(\TT^{n+1}\right)$, that is
	\[
	\sum_{\left(\tau,\xi\right)\in\Z^{n+1}} \left|\hat u(\tau,\xi)\right|^2\left(1+\left|(\tau,\xi)\right|\right)^{-2p} <\infty.
	\]
	Since $Pu\in\cinfty\left(\TT^{n+1}\right)$, from Corollary \ref{Cor:ellipticity-wf} we obtain that $\hat u(\cdot,\xi) \in \cinfty\left(\TT^n\right)$ for every $\xi \in \Z$ and there is $c>0$ such that given $k\in\Z_+$ there exists $C_k>0$ such that
	\[
	\left|\hat u(\tau,\xi)\right| \le C_k\left(1+|(\tau,\xi)|\right)^{-k}, \quad \forall|\xi|\le c|\tau|.
	\] 
	Then 
	\begin{align*}
		\left|\hat u(t,\xi)\right| &= \left|\sum_{\tau \in \Z^n} \hat u(\tau,\xi) e^{it\tau}\right| \le \sum_{\tau \in \Z^n} \left|\hat u(\tau,\xi)\right|\\
		&=\sum_{\tau \in \Z^n} \left(1+|\tau|\right)^{-n}\left(1+|\tau|\right)^{n}\left|\hat u(\tau,\xi)\right|\\
		&\le \left(\sum_{\tau \in \Z^n}\left(1+|\tau|\right)^{-2n}\right)^{1/2} \left(\sum_{\tau\in\Z^n} \left(1+|\tau|\right)^{2n}\left|\hat u(\tau,\xi)\right|^2\right)^{1/2}.
	\end{align*}
	Now we observe that
	\begin{align*}
		\sum_{\tau\in\Z^n} \left(1+|\tau|\right)^{2n}\left|\hat u(\tau,\xi)\right|^2 &= \sum_{|\xi|\le c|\tau|} \left(1+|\tau|\right)^{2n}\left|\hat u(\tau,\xi)\right|^2+\sum_{|\xi|>c|\tau|} \left(1+|\tau|\right)^{2n}\left|\hat u(\tau,\xi)\right|^2\\
		&= (I)+(II)
	\end{align*}
	and we have:
	\begin{align*}
		(I) &\le \sum_{|\xi|\le c|\tau|} \left(1+|\tau|\right)^{2n} C_{2n} \left(1+|(\tau,\xi)|\right)^{-4n}\\
		&\le  C_{2n}\sum_{\tau \in \Z^n} \left(1+|\tau|\right)^{-2n}\\
		&\le  \left(C_{2n}\sum_{\tau \in \Z^n} \left(1+|\tau|\right)^{-2n}\right)\left(1+|\xi|\right)^{2(n+p)}
	\end{align*}
	and, for some $c_1>0$,
	\begin{align*}
		(II) &= \sum_{|\xi|>c|\tau|} \left(1+|\tau|\right)^{2n}\left(1+|(\tau,\xi)|\right)^{2p-2p}\left|\hat u(\tau,\xi)\right|^2\\
		&\le c_1\left(1+|\xi|\right)^{2(n+p)}\sum_{\tau\in\Z^n}\left(1+|(\tau,\xi)|\right)^{-2p}\left|\hat u(\tau,\xi)\right|^2
	\end{align*}
	from which we can conclude that there exists $D>0$ such that
	\[
	\left|\hat u(t,\xi)\right| \le D\left(1+|\xi|\right)^{2(n+p)}, \quad \forall t \in \TT^n, \xi \in \Z.
	\]
\end{proof}
\begin{Lem}
	\label{Lem:exp-CD-estimate}
	Let $\alpha \in \R, \ell \in \Z_+$ and $\xi \in \Z$ such that
	\[
	\left|q-\alpha\xi\right| \ge \left(1+|\xi|\right)^{-\ell}, \quad \forall q \in \Z.
	\]
	Then
	\[
	\left|e^{2\pi i \alpha\xi}-1\right| \ge \left(1+|\xi|\right)^{-\ell}.
	\]
\end{Lem}
\begin{proof}
	One of the following conditions occurs:
	\begin{enumerate}
		\item There exists $q\in\Z$ such that $\left|\pi \alpha \xi - \pi q\right| < \frac{\pi}{3}$;
		\item For every $q\in\Z$, we have $\left|\pi\alpha\xi - \pi q\right| \ge \frac{\pi}{3}$.
	\end{enumerate}
	Suppose that the first condition holds true. Then for some $\theta$ between $\pi\alpha\xi$ and $\pi q$ we have
	\begin{align*}
		\left|e^{2\pi i \alpha\xi}-1\right| &= 2\left|\sin\left(\pi\alpha\xi\right)\right|\\
		&= 2\left|\sin\left(\pi\alpha\xi\right)-\sin\left(q\pi\right)\right|\\
		&=2\left|\cos(\theta)\right| \left|\pi\alpha\xi-q\pi\right|\\
		&\ge \pi\left|\alpha\xi-q\right|\\
		&\ge \left(1+|\xi|\right)^{-\ell}
	\end{align*}
	since $|\cos\left(\theta\right)| \ge \left|\cos\left(\frac{\pi}{3}\right)\right| = \frac{1}{2}$.
	
	If the second condition occurs then
	\begin{align*}
		\left|e^{2\pi i \alpha\xi}-1\right| &= 2\left|\sin\left(\pi\alpha\xi\right)\right| \ge \sqrt{3} \ge	\left(1+|\xi|\right)^{-\ell}.
	\end{align*}
\end{proof}
Now we can prove the central point of this work.
%
%
\\

\noindent{\bf Proof of Theorem \ref{Thm:propagation-one-point}}: By Lemma \ref{Lem:sobolev-estimate} there are $M \in \Z_+$ and $D>0$ such that
	\begin{align}
		\label{eq:thm-prop-one-point-sobolev}
		\left|\widehat u(t,\xi)\right| \le D\left(1+|\xi|\right)^{M},\quad \forall t \in \TT^n, \xi \in \Z.
	\end{align}
	If we set $f_j=X_ju, j \in\{1,\ldots,n\}$, then given $k\in\Z_+$ there exists $D_k>0$ such that
	\[
	\left|\widehat{S_j f_j}(t,\xi)\right| = \left|\widehat{f_j}(t,\xi)\right| \le D_k\left(1+|\xi|\right)^{-2k-M}, \quad \forall t\in \TT^n, j \in\{1,\ldots,n\}.
	\]
	We claim that if $B_k=\max\{C_k,D_k\}$ and $C_k' = 2\pi (n+1)\left(B_k+D_k+D\right)$ then for every $t \in \TT^n$, $k\in\Z_+$ and $\xi \in \Z$ we have
	\begin{align}
		\label{eq:thm-prop-oone-point-goal}
		\left|\hat u(t,\xi)\right| \le C_k'\left(1+|\xi|\right)^{-k}.
	\end{align}
	%
	%
	%

	%
	To prove this inequality we fix $\xi \in \Z, k\in\Z_+$ and we will show that every $t \in \TT^n$ satisfies \eqref{eq:thm-prop-oone-point-goal}. Note that if $t \in \TT^n$ is such that there is $j \in \{1,\ldots,n\}$ satisfying
	\begin{align}
		\label{eq:set-A-xi}
		\left|q - a_{j0}(t)\xi\right| \ge \left(1+|\xi|\right)^{-k-M}, \quad \forall q\in\Z,
	\end{align}
	then 
	\begin{align*}
		\left|\hat u(t,\xi)\right| \le 2\pi B_k\left(1+|\xi|\right)^{-k}
	\end{align*}
	and in particular \eqref{eq:thm-prop-oone-point-goal} holds true for $t$. Indeed, by Lemma \ref{Lem:exp-CD-estimate} we obtain that
	\[
	\left|e^{2\pi i a_{j0}(t)\xi}-1\right| \ge \left(1+|\xi|\right)^{-k-M}
	\]
	and we use Lemma 3.2 of~\cite{hps06} (recall \eqref{eq:automorphism_Sj-fourier-formula} and that $a_{j0}$ does not depend on $t_j$) to write
	\[
	\widehat{S_ju}\left(t,\xi\right) = \frac{1}{e^{2\pi i a_{j0}(t)\xi}-1}\int_0^{2\pi} e^{i s a_{j0}(t)\xi} \widehat{S_j f_j}\left(t_1,\ldots,t_{j-1},t_j+s,t_{j+1},\ldots,t_n,\xi\right)\dd s,
	\]
	which yields
	\begin{align*}
		\left|\hat u(t,\xi)\right| &= \left|\widehat{S_ju}\left(t,\xi\right)\right|\\
		&\le \frac{1}{\left(1+|\xi|\right)^{-k-M}} 2\pi \sup_{t' \in \TT^n} \left|\widehat{S_jf_j}\left(t',\xi\right)\right|\\
		&=  \frac{2\pi}{\left(1+|\xi|\right)^{-k-M}}  \sup_{t' \in \TT^n} \left|\widehat{f_j}\left(t',\xi\right)\right|\\
		&\le 2\pi B_k\left(1+|\xi|\right)^{-k}.
	\end{align*}
	Now define $\mathcal A$ as the set of the points $t \in \TT^n$ such that there exists $j \in \{1,\ldots,n\}$ satisfying \eqref{eq:set-A-xi} and $\mathcal B \doteq \TT^n\setminus \mathcal A$. Hence $t \in \mathcal B$ if and only if given $j \in \{1,\ldots,n\}$, there is $q \in \Z$ such that
	\[
	\left|q-a_{j0}(t)\xi\right| < \left(1+|\xi|\right)^{-k-M}.
	\]
	We then consider, for each $\tau \in \Z^n$, the set $\mathcal B_{\tau}$ of the points $t \in \TT^n$ such that
	\[
	\max_{1\le j\le n} \left|\tau_j-a_{j0}(t)\xi\right| < \left(1+|\xi|\right)^{-k-M}
	\]
	and we have that $\mathcal B = \cup_{\tau \in \Z^n}\mathcal B_{\tau}$. Since every $t \in \mathcal A$ satisfies \eqref{eq:thm-prop-oone-point-goal}, to finish the proof we will fix $\tau \in \Z^n$ and show that \eqref{eq:thm-prop-oone-point-goal} holds true for each point in $\mathcal B_{\tau}$. 
	%
	Fix $\tau \in \Z^n$ and $t = \left(t_1,\ldots,t_n\right) \in \mathcal B_\tau$. Recall that $a_{j0}$ does not depend on $t_j$ so we can consider
	\[
	I_n = \left\{t_n' \in \TT^{n-1} : \left|\tau_n - a_{n0}\left(t_n'\right)\xi\right|<\left(1+|\xi|\right)^{-k-M}\right\},
	\]
	which contains $\left(t_1,\ldots,t_{n-1}\right)$ and, since for any $\varphi \in \cinfty\left(\TT^n\right)$
	\begin{align*}
		\varphi\left(t\right) &= \varphi\left(t_1,\ldots,t_{n-1},s_n\right) + \int_{s_n}^{t_n}\partial_n \varphi\left(t_1,\ldots,t_{n-1},\lambda\right)\dd\lambda,
	\end{align*}
	we have that
	\[
	\left|\varphi\left(t\right)\right| \le \left|\varphi\left(t_1,\ldots,t_{n-1},s_n\right)\right| + 2\pi \sup_{(t'_n,\lambda)\in I_n\times \TT}\left|\partial_n\varphi\left(t'_n,\lambda\right)\right|.
	\]
	We use the last inequality with 
	\[
	\varphi(y) = e^{iy_n \tau_n}\widehat{S_n u}\left(y,\xi\right)
	\]
	 and it follows from \eqref{eq:partial-fourier-auto-Sj-modulos} that
	\[
	\left|\hat u\left(t,\xi\right)\right| \le \left|\hat u\left(t_1,\ldots,t_{n-1},s_n,\xi\right)\right| + 2\pi \sup_{(t'_n,\lambda)\in I_n\times \TT}\left|\left(\partial_{t_n} + i\tau_n\right)\widehat{S_n u}(t'_n,\lambda,\xi)\right|
	\]
	By using \eqref{eq:thm-prop-one-point-sobolev} and \eqref{eq:automorphism_Sj-fourier-formula} we conclude, for $t'_n \in I_n$ and $\lambda \in \TT$, that
	\begin{align*}
		\left|\left(\partial_{t_n} + i\tau_n\right)\widehat{S_n u}(t'_n,\lambda,\xi)\right|  &\le \left|\left(\partial_{t_n} + ia_{n0}\left(t'_n\right)\xi\right)\widehat{S_n u}(t'_n,\lambda,\xi)\right|\\
		&+\left|\tau_n-a_{n0}\left(t'_n\right)\xi\right| \left|\widehat{S_n u}\left(t_n',\lambda,\xi\right)\right|\\
		&\le \left|\widehat{S_nf_n}\left(t'_n,\lambda,\xi\right)\right| + D\left(1+|\xi|\right)^{-k}\\
		&\le \left(D_k + D\right)\left(1+|\xi|\right)^{-k}.
	\end{align*}
	Now if $\left(t_1,\ldots,t_{n-1},s_n\right)\in \mathcal A$ then
	\[
	\left|\hat u\left(t,\xi\right)\right| \le 2\pi\left(B_k+D_k + D\right)\left(1+|\xi|\right)^{-k}
	\]
	and we obtain \eqref{eq:thm-prop-oone-point-goal}. Otherwise, there exists $\tilde\tau = \left(\tilde \tau_1,\ldots,\tilde \tau_n\right)\in\Z^n$ such that $\left(t_1,\ldots,t_{n-1},s_n\right)\in B_{\tilde\tau}$ and we consider
	\[
	I_{n-1} = \left\{t'_{n-1} \in \TT^{n-1}: \left|\tilde\tau_{n-1} - a_{(n-1)0}\left(t'_{n-1}\right)\xi\right|<\left(1+|\xi|\right)^{-k-M}\right\},
	\]
	which contains $\left(t_1,\ldots,t_{n-2},s_n\right)$. Since for any $\varphi \in \cinfty\left(\TT^n\right)$
	\begin{align*}
		\varphi\left(t_1,\ldots,t_{n-1},s_n\right) &= \varphi\left(t_1,\ldots,t_{n-2},s_{n-1},s_n\right) + \int_{s_{n-1}}^{t_{n-1}}\partial_{n-1} \varphi\left(t_1,\ldots,t_{n-2},\lambda,s_n\right)\dd\lambda,
	\end{align*}
	we have that
	\begin{align*}
	\left|\varphi\left(t_1,\ldots,t_{n-1},s_n\right)\right| &\le \left|\varphi\left(t_1,\ldots,t_{n-2},s_{n-1},s_n\right)\right|+ \\
	&\quad+ 2\pi \sup_{(t'_{n-1},\lambda)\in I_{n-1}\times \TT}\left|\partial_{n-1}\varphi\left(t_1,\ldots,t_{n-2},\lambda,t_{n}\right)\right|.
	\end{align*}
	We use this inequality with 
	\[
	\varphi(y) = e^{iy_{n-1}\tilde \tau_{n-1}}\widehat{S_{n-1}u}(y,\xi)
	\] 
	to obtain
	\begin{align*}
		\left|\hat u\left(t_1,\ldots,t_{n-1},s_n,\xi\right)\right| &\le \left|\hat u\left(t_1,\ldots,t_{n-2},s_{n-1},s_n,\xi\right)\right|+\\ 	
		& + 2\pi \sup_{(t'_{n-1},\lambda)\in I_{n-1}\times \TT}\left|\left(\partial_{n-1}+i\tilde \tau_{n-1}\right)\widehat{S_{n-1}u}\left(t_1,\ldots,t_{n-2},\lambda,t_{n},\xi\right)\right|.
	\end{align*}
	As before, we conclude that
	\[
	\sup_{\left(t_{n-1}',\lambda\right) \in I_{n-1}\times \TT}\left|\left(\partial_{n-1}+i\tilde \tau_{n-1}\right)\widehat{S_{n-1}u}\left(t_1,\ldots,t_{n-2},\lambda,t_{n},\xi\right)\right| \le \left(D_k+D\right)\left(1+|\xi|\right)^{-k},
	\]
	%
so
	\begin{align*}
		\left|\hat u\left(t,\xi\right)\right| &\le \left|\hat u\left(t_1,\ldots,t_{n-1},s_n,\xi\right)\right| + 2\pi \left(D_k + D\right)\left(1+|\xi|\right)^{-N}\\
		&\le \left|\hat u\left(t_1,\ldots,t_{n-2},s_{n-1},s_n,\xi\right)\right| + 4\pi \left(D_k + D\right)\left(1+|\xi|\right)^{-k}.
	\end{align*}
	If $\left(t_1,\ldots,t_{n-2},s_{n-1},s_n\right) \in \mathcal A$, we obtain~\eqref{eq:thm-prop-oone-point-goal}. Otherwise we proceed with the same argument: either there is $j\in\{1,\ldots,n\}$ such that  $\left(t_1,\ldots,t_j,s_{j+1},\ldots,s_n\right) \in \mathcal A$ or we use the hypothesis about $s$ to conclude the proof.	
\begin{flushright}
	\qed
\end{flushright}

\section{A connection between sum of squares and systems}
\label{sec:system-sums}

In this section we shall prove that the system $\mathcal X$ defined by \eqref{eq:intro-system} is globally hypoelliptic if and only if its sum of squares associated $P=\sum_{j=1}^n X_j^2$ is globally hypoelliptic. We recall that a system $\mathcal S$ of linear partial differential operators in $\TT^N$ is globally hypoelliptic in $\TT^N$ when $u \in \D'\left(\TT^N\right)$ and $Qu\in\cinfty\left(\TT^N\right)$ for every $Q \in \mathcal S$ ensure that $u\in\cinfty\left(\TT^N\right)$.
\\

\noindent{\bf Proof of Theorem~\ref{Thm:P-GH-iff-L-GH}}: It is clear that if $u \in \D'\left(\TT^{n+1}\right)$ and $X_j u \in \cinfty\left(\TT^{n+1}\right)$ for every $j \in \{1,\ldots,n\}$ then $Pu \in \cinfty\left(\TT^{n+1}\right)$. Conversely, suppose that $Pu=f \in \cinfty\left(\TT^{n+1}\right)$ for some $u \in \D'\left(\TT^{n+1}\right)$ and let us prove that $X_ju \in \cinfty\left(\TT^{n+1}\right)$ for every $j\in\{1,\ldots,n\}$. Taking the Fourier transform with respect to $x$, we obtain
	\[
	\sum_{j=1}^n \left(\partial_{t_j} + i a_j(t)\xi\right)^2\hat u(t,\xi) = \hat f(t,\xi), \quad \forall \xi \in \Z.
	\]
	By ellipticity we conclude that $\hat u(\cdot,\xi) \in \cinfty\left(\TT^{n}\right)$ for every $\xi \in \Z$. Moreover, by Corollary \ref{Cor:ellipticity-wf} there exists $c>0$ such that for any $k\in\Z_+$, one can find $C_k>0$ satisfying
	\[
	\left|\hat u(\tau,\xi)\right| \le C_k\left(1+|(\tau,\xi)|\right)^{-k},\quad \forall |\xi| \le c|\tau|.
	\]
	By using Proposition \ref{Prop:inclusion_su} we can increase $C_k$ and decrease $c$ if necessary to obtain
	\[
	\left|\widehat {X_ju}(\tau,\xi)\right| \le C_k\left(1+|(\tau,\xi)|\right)^{-k},\quad \forall |\xi| \le c|\tau|
	\]
	for every $j \in \{1,\ldots,n\}$. It remains to prove a similar inequality for frequencies satisfying $|\xi| > c |\tau|$. Since $|\widehat{X_j u}(\tau,\xi)| \le \left\|\widehat{X_ju}(\cdot,\xi)\right\|_{L^2\left(\TT^n\right)}$, it suffices to prove that for any $k\in\Z_+$ there is $C_k'>0$ such that
	\[
	\left\|\widehat{X_ju}(\cdot,\xi)\right\|_{L^2\left(\TT^n\right)}^2 \le C_k'\left(1+|\xi|\right)^{-k}, \quad \forall \xi \in \Z.
	\]
	This can be done in the following way: if we set
	\[
	X_{j\xi} \doteq \partial_{t_j} + ia_j(t)\xi, \quad \xi\in\Z,
	\]
	then
	\begin{align*}
		\sum_{j=1}^n \left\|\widehat{X_ju}(\cdot,\xi)\right\|_{L^2\left(\TT^n\right)}^2 &= \left|\sum_{j=1}^n \left\langle \widehat{X_ju}(\cdot,\xi),\widehat{X_ju}(\cdot,\xi) \right\rangle_{L^2\left(\TT^n\right)}\right|\\
		&= \left|\sum_{j=1}^n \left\langle X_{j\xi}\hat u(\cdot,\xi),X_{j\xi}\hat u(\cdot,\xi) \right\rangle_{L^2\left(\TT^n\right)}\right|\\
		&=  \left|\sum_{j=1}^n -\left\langle X_{j\xi}^2\hat u(\cdot,\xi),\hat u(\cdot,\xi) \right\rangle_{L^2\left(\TT^n\right)}\right|\\
		&= \left| \left\langle  \widehat{Pu}(\cdot,\xi),\hat u(\cdot,\xi) \right\rangle_{L^2\left(\TT^n\right)}\right|\\
		&\le \left\| \widehat{Pu}(\cdot,\xi)\right\|_{L^2\left(\TT^n\right)} \left\|u(\cdot,\xi)\right\|_{L^2\left(\TT^n\right)}.
	\end{align*}
	We know that for each $k\in\Z_+$ there is $D_k>0$ such that 
 \[
 \left\| \widehat{Pu}(\cdot,\xi)\right\|_{L^2\left(\TT^n\right)} \le D_k\left(1+|\xi|\right)^{-k},\quad \forall\xi \in \Z.
 \]
 Thus it suffices to prove that there are $C,M>0$ such that
	\[
	\left\|\hat u\left(\cdot,\xi\right)\right\|_{L^2\left(\TT^n\right)} \le C\left(1+|\xi|\right)^M, \quad \forall \xi \in \Z,
	\]
	which follows from Lemma~\ref{Lem:sobolev-estimate}.
	\begin{flushright}
		\qed
	\end{flushright}
%
%
%


%
%
\begin{Rem}
	\label{Rem:petronilho-conj}
	If we consider, in $\TT_t^n\times\TT_x^m$,
	\[
	P=-\sum_{j=1}^n\left(\partial_{t_j} + \sum_{k=1}^m a_{jk}(t)\partial_{x_k}\right)^2,
	\]
	where $a_{jk} \in \cinfty\left(\TT^n\right)$ are real, then we can adapt the proof of Theorem \ref{Thm:P-GH-iff-L-GH} to show that $P$ is globally hypoelliptic if and only if the system $\mathcal X=\left\{\partial_{t_j} + \sum_{k=1}^m a_{jk}(t)\partial_{x_k}\right\}_{j=1}^n$ is globally hypoelliptic. In particular, the global hypoellipticity of one of the vector fields of $\mathcal X$ ensures the global hypoellipticity of $P$. Therefore, Theorem 3.3 and Corollary 3.5 of \cite{p06} are an immediate consequence of Theorem \ref{Thm:P-GH-iff-L-GH}. Also, for this class of operators, Petronilho's conjecture \cite{p06} has an affirmative answer.
\end{Rem}
\begin{Rem}
	\label{Rem:sums-squares-with-laplacian}
	In $\TT^n_t\times\TT^m_x$ we consider
	\[
	P = -\sum_{j=1}^n \partial_{t_j}^2 - \sum_{\ell=1}^N\left(\sum_{k=1}^m a_{\ell k}(t)\partial_{x_k} + W_\ell\right)^2,
	\]
	where $W_\ell$ is a skew-symmetric real vector field in $\TT^n$ and $N\in\N$ is arbitrary. We denote 
	\[
	X_\ell = \sum_{k=1}^m a_{\ell k}(t)\partial_{x_k} + W_\ell, \quad \ell \in \{1,\ldots,N\},
	\]
	and one can adapt the argument of Theorem \ref{Thm:P-GH-iff-L-GH} in order to show that $P$ is globally hypoelliptic in $\TT^{n+m}$ if and only if
	\[
	\mathcal L = \left\{\partial_{t_1},\ldots,\partial_{t_n},X_1,\ldots,X_n\right\}
	\]
	is globally hypoelliptic in $\TT^{n+m}$. If we suppose that $Xu \in \cinfty\left(\TT^{n+m}\right)$ for every $X \in \mathcal L$, then $\partial_{t_j} u  \in \cinfty\left(\TT^{n+m}\right)$, so $W_j u \in \cinfty\left(\TT^{n+m}\right)$. In particular, we obtain that
	\[
	\sum_{k=1}^m a_{\ell k}(t)\partial_{x_k}u \in \cinfty\left(\TT^{n+m}\right), \quad \forall j \in \{1,\ldots,n\}.
	\]
	Also, there is a regularity in $t$-variable (in the sense of Corollary~\ref{Cor:ellipticity-wf}). Therefore it is natural that the global hypoellipticity of $P$ in $\TT^{n+m}$ is equivalent to a condition about the functions $a_{\ell k}$ in $\TT^n$ as it occurs in the previous works (see \cite{bfp17,hp00} for instance).
\end{Rem}

\section{Global Hypoellipticity}
\label{sec:GH}

To state our result about global hypoellipticity we need the following diophantine condition (see~\cite{bcm93,hp99}): we say that a collection of real numbers $\{\alpha_1,\ldots,\alpha_n\} \subset \R^m$ is not simultaneously approximable when there exist $C,\rho>0$ such that
\[
\max_{1\le j \le n} \left|\eta_j + \alpha_j \xi\right| \ge C\left(1+|\xi|\right)^{-\rho}, \quad \forall \eta=\left(\eta_1,\ldots,\eta_n\right) \in \Z^n, \xi \in \Z\setminus\{0\}.
\]
Otherwise we say that $\left\{\alpha_1,\ldots,\alpha_n\right\}$ is simultaneously approximable.
\begin{Prop}
	\label{Prop:SA}
	A collection of real numbers $\{\alpha_1,\ldots,\alpha_n\}$ is simultaneously approximable if and only if there exists a sequence $\left(\tau^\nu,\xi^\nu\right) \in \Z^n\times\Z$ such that $\xi^\nu\ne 0$ for every $\nu \in \N$, $\left|(\tau^\nu,\xi^\nu)\right| \to \infty$ and
	\[
	\max_{1\le j\le n} \left|\tau_j^\nu +\alpha_j\xi^\nu\right| < \frac{C_0}{\left(1+\left|(\tau^\nu,\xi^\nu)\right|\right)^\nu}, \quad \forall \nu \in \N
	\]
	for some $C_0>0$.
\end{Prop}
\begin{proof}
	Suppose that $\left\{\alpha_1,\ldots,\alpha_n\right\}$ is simultaneously approximable. We split the proof into cases:
	\begin{enumerate}
		\item There exists $\left(\tau,\xi\right) \in \Z^n\times \Z$ such that $\xi\ne0$ and
		\[
		\tau_j +\alpha_j\xi = 0, \quad \forall j \in \{1,\ldots,n\}.
		\]
		In this case we can take $\tau^\nu = \nu \tau$ and $\xi^\nu=\nu \xi$ for every $\nu \in \N$.
		
		\item For every $\left(\tau,\xi\right)\in\Z^n\times \Z$ with $\xi\ne 0$ we have that $\max_{1\le j\le n} \left|\tau_j +\alpha_j\xi\right| \ne 0$. Then given $C=1$ and $\rho=2$, there exists $\left(\tau^1,\xi^1\right) \in \Z^n\times \Z, \xi^1\ne 0$, such that
		\[
		0<\max_{1\le j\le n} \left|\tau_j^1 +\alpha_j\xi^1\right| < \frac{1}{\left(1+\left|\xi^1\right|\right)^2}.
		\]
		Now we take 
		\[
		0<C<\min\left\{1,\max_{1\le j\le n} \left|\tau_j^1 +\alpha_j\xi^1\right|\left(1+\left|\xi^1\right|\right)^2\right\}
		\]
		and $\rho=4$, so there exists $\left(\tau^2,\xi^2\right)\in\Z^n\times\Z$ with $\xi^2\ne 0$ and
		\[
		0<\max_{1\le j\le n} \left|\tau_j^2 +\alpha_j\xi^2\right| < \frac{C}{\left(1+\left|\xi^2\right|\right)^4} \le \frac{C}{\left(1+\left|\xi^2\right|\right)^2}.
		\]
		By the choice of $C$ we obtain that $\left(\tau^1,\xi^1\right) \ne\left(\tau^2,\xi^2\right)$ and that
		\[
		\max_{1\le j\le n} \left|\tau_j^2 +\alpha_j\xi^2\right| < \frac{1}{\left(1+\left|\xi^2\right|\right)^4}.
		\]
		Proceeding in this way we obtain a sequence $\left(\tau^\nu,\xi^\nu\right)\in\Z^n\times\Z$ such that $\xi^\nu\ne 0$ for every $\nu \in \N$, $\left(\tau^\nu,\xi^\nu\right)\ne \left(\tau^{\nu'},\xi^{\nu'}\right)$ if $\nu\ne \nu'$ and
		\[
		\max_{1\le j\le n} \left|\tau_j^\nu +\alpha_j\xi^\nu\right| < \frac{1}{\left(1+\left|\xi^\nu\right|\right)^{2\nu}}, \quad \forall \nu \in \N.
		\]
		In particular, $\left|\left(\tau^\nu,\xi^\nu\right)\right|\to \infty$ and if $\alpha=\max_{1\le j\le n} \left|\alpha_j\right|$ then
		\begin{align*}
			\left|\tau^\nu_j\right| &\le \left|\tau^\nu_j +\alpha_j \xi^\nu\right| + \alpha\left|\xi^\nu\right|\\
			&\le \frac{1}{\left(1+\left|\xi^\nu\right|\right)^\nu} + \alpha\left|\xi^\nu\right|\\
			&\le 1 + \alpha\left|\xi^\nu\right|\\
			&\le \left(1+\alpha\right) \left|\xi^\nu\right|, \quad \forall j \in \{1,\ldots,n\}.
		\end{align*}
		Thus
		\begin{align*}
			\left|\left(\tau^\nu,\xi^\nu\right)\right|^2 &= \sum_{j=1}^n \left|\tau_j^\nu\right|^2 + \left|\xi^\nu\right|^2\\
			&\le \sum_{j=1}^n \left(1+\alpha\right)^2\left|\xi^\nu\right|^2 + \left|\xi^\nu\right|^2\\
			&= \left[n(1+\alpha)^2+1\right]\left|\xi^\nu\right|^2.
		\end{align*}
		If $c_0 = \sqrt{n(1+\alpha)^2+1}$ then
		\[
		1+\left|\left(\tau^\nu,\xi^\nu\right)\right| \le \left(1+c_0\right)\left|\xi^\nu\right| \le \left(1+c_0\right)\left(1+\left|\xi^\nu\right|\right),
		\]
		which yields
		\[
		\frac{1}{\left(1+\left|\xi^\nu\right|\right)^{2\nu}} \le \frac{1}{\left(1+\left|\left(\tau^\nu,\xi^\nu\right)\right|\right)^\nu}\frac{\left(1+c_0\right)^{2\nu}}{\left(1+\left|\left(\tau^\nu,\xi^\nu\right)\right|\right)^\nu}, \quad \forall \nu \in \N.
		\]
		Since $\left|\left(\tau^\nu,\xi^\nu\right)\right|\to \infty$, there exists $\nu_0$ such that
		\[
		\left(1+c_0\right)^2 \le \left|\left(\tau^\nu,\xi^\nu\right)\right|, \quad \forall \nu \ge \nu_0
		\]
		and starting the sequence from $\nu_0$ we obtain the desired sequence.
	\end{enumerate}
	
	Conversely, suppose the existence of such sequence. 
	Thus given $C,\rho>0$ there exists $\nu \in \N$ such that
	\[
	\max_{1\le j \le n} \left|\tau_j^\nu - \alpha_j\xi^\nu\right| \le \frac{C_0}{\left(1+\left|\left(\tau^\nu,\xi^\nu\right)\right|\right)^\nu} \le \frac{C}{\left(1+\left|\left(\tau^\nu,\xi^\nu\right)\right|\right)^\rho}
	\]
	and the proof is complete.

\end{proof}
\noindent{\bf Proof of Theorem~\ref{Thm:equivalence-gh}}: suppose that there is no point $t \in \TT^n$ of finite type for $X_1,\ldots,X_n$ and that the range of each $a_{10},\ldots,a_{n0}$ does not contain real numbers $\alpha_1,\ldots,\alpha_n \in \R$, respectively, that are not simultaneously approximable. Our goal is to prove that $\mathcal X$ is not globally hypoelliptic in $\TT^{n+1}$. By Remark \ref{Rem:infinite-type} we know that there are $\alpha_1,\ldots,\alpha_n \in \R$ such that $a_{j0}=\alpha_j$ for every $t \in \TT^n$. By hypothesis, $\alpha_1,\ldots,\alpha_n \in \R$ are simultaneously approximable: there are $C_0>0$ and a sequence $\left(\tau^\nu,\xi^\nu\right) \in \Z^{n+1}$ such that $\left|\left(\tau^\nu,\xi^\nu\right)\right| \to \infty$ and
\[
\max_{1\le j\le n} \left|\tau^\nu_j + \alpha_j \xi^\nu\right| < \frac{C_0}{\left(1+\left|\left(\eta^\nu,\xi^\nu\right)\right|\right)^\nu}, \quad \forall \nu \in \N.
\]
We set
\[
v \doteq \sum_{\nu=1}^\infty e^{i(t,x)\cdot\left(\tau^\nu,\xi^\nu\right)} \in \D'\left(\TT^{n+1}\right)\setminus \cinfty\left(\TT^{n+1}\right).
\]
It follows from Remark \ref{Rem:infinite-type} that there is an automorphism $S:\D'\left(\TT^{n+1}\right)\to \D'\left(\TT^{n+1}\right)$ that restricts itself to an automorhpism of $\cinfty\left(\TT^{n+1}\right)$ that satisfies $Y_j S = SX_j$ for every $j \in\{1,\ldots,n\}$, where
\[
Y_j = \partial_{t_j} + \alpha_j\partial_x.
\]
We set $u=S^{-1}v \in\D'\left(\TT^{n+1}\right) \setminus \cinfty\left(\TT^{n+1}\right)$ and claim that $X_j u \in \cinfty\left(\TT^{n+1}\right)$ for every $j\in\{1,\ldots,n\}$. Indeed,
\begin{align*}
	X_j u &= X_j S^{-1} v = S^{-1} Y_j v\\
	&= S^{-1} \left(i\sum_{\nu=1}^\infty \left(\tau_j^\nu + \alpha_j\xi^\nu\right)e^{i(t,x)\cdot\left(\tau^\nu,\xi^\nu\right)}\right) \in \cinfty\left(\TT^{n+1}\right).
\end{align*}
Therefore $\mathcal X$ is not globally hypoelliptic in $\TT^{n+1}$.

Conversely, let $u\in\D'\left(\TT^{n+1}\right)$ such that $X_j u \in \cinfty\left(\TT^{n+1}\right)$ for every $j\in\{1,\ldots,n\}$ and we shall use Theorem \ref{Thm:propagation-one-point} in order to prove that $u \in \cinfty\left(\TT^{n+1}\right)$ in the case that there exists $t \in \TT^n$ of finite type for $X_1,\ldots,X_n$. In this case there are $j,\ell \in \{1,\ldots,n\}$ such that
\[
\left[X_j,X_\ell\right] = \left(\partial_{t_j}a_{\ell}-\partial_{t_\ell}a_{j}\right)\partial_{x} \ne 0.
\]
Hence one can find an open subset $U\subset \TT^n$ such that
\[
b(t) \doteq \frac{1}{\partial_{t_j}a_{\ell}(t)-\partial_{t_\ell}a_{j}(t)}  \in \cinfty\left(U\right)\cap L^\infty\left(U\right).
\]
Since $f\doteq \left[X_j,X_\ell\right] u \in \cinfty\left(\TT^{n+1}\right)$, given $k\in\Z_+$, there exists $C_k>0$ such that
\[
\left|\hat f(t,\xi)\right| \le C_k\left(1+|\xi|\right)^{-k}, \quad \forall t \in \TT^n, \xi \in \Z.
\]
Then for $t \in U$ and $\xi \in \Z\setminus\{0\}$,
\begin{align*}
	\left|\hat u(t,\xi)\right| &= \left|\frac{\xi}{\xi}\hat u(t,\xi)\right|\\
	&= \left|\frac{1}{\xi}\widehat{\partial_x u}(t,\xi)\right|\\
	&= \frac{1}{|\xi|}\left|b(t)\hat{f}(t,\xi)\right|\\
	&\le \left\|b\right\|_{L^\infty(U)} C_{k+1}\left(1+|\xi|\right)^{-k}
\end{align*}
and by Theorem \ref{Thm:propagation-one-point} we conclude that $u \in \cinfty\left(\TT^{n+1}\right)$. 

Now suppose that there is no point of finite type for $X_1,\ldots,X_n$ and the range of each $a_{10},\ldots,a_{n0}$ contains real numbers $\alpha_1,\ldots,\alpha_n \in \R^m$, respectively, that are not simultaneously approximable: there are $C,\rho>0$ such that
\[
\max_{1\le j \le n} \left|\tau_j + \alpha_j\xi\right| \ge C\left(1+\left|\xi\right|\right)^{-\rho}, \quad \forall \left(\tau,\xi\right)\in\Z^n\times\Z, \xi \ne0.
\]
Recall (see Remark~\ref{Rem:infinite-type}) that in this case, $a_{j0}(t) = \alpha_j$ for every $j \in\{1,\ldots,n\}$ and $t \in \TT^n$ and we also consider $S$ that satisfies \eqref{eq:simultaneous-reduction}. It suffices to prove that $v\doteq Su \in \cinfty\left(\TT^{n+1}\right)$. Since
\[
g_j \doteq Y_j S u = S X_j u \in \cinfty\left(\TT^{n+1}\right), \quad \forall j \in\{1,\ldots,n\},
\] 
given $k\in\Z_+$ there exists $C_k>0$ such that
\[
\left|\widehat{g_j}(\tau,\xi)\right| \le C_k\left(1+|(\tau,\xi)|\right)^{-k}, \quad \forall (\tau,\xi) \in \Z^{n+1}, j \in\{1,\ldots,n\}.
\]
Since
\[
i\left(\tau_j+\alpha_j\xi\right) \hat v(\tau,\xi) = \widehat{g_j}(\tau,\xi), \quad \forall j \in\{1,\ldots,n\}, \left(\tau,\xi\right)\in\Z^{n+1},
\]
we obtain that
\[
\left|\hat v(\tau,\xi)\right| \le C^{-1}\left(1+|\xi|\right)^\rho C_k\left(1+|(\tau,\xi)|\right)^{-k}, \quad \forall \tau \in \Z^n, \xi \ne 0.
\]
It remains to prove a similar inequality when $\xi=0$. Since $Pu \in \cinfty\left(\TT^{n+1}\right)$ and $Y_j S = S X_j$ for every $j\in\{1,\ldots,n\}$, it follows that
\[
\sum_{j=1}^n Y_j^2 v = S P u \in \cinfty\left(\TT^{n+1}\right)
\]
and from Corollary \ref{Cor:ellipticity-wf} (for $\sum_{j=1}^n Y_j^2$) we conclude that, for some $c>0$, given $k\in\Z_+$ there exists $D_k>0$ such that
\[
\left|\hat v(\tau,\xi)\right| \le D_k\left(1+|(\tau,\xi)|\right)^{-k}, \quad \forall |\xi|\le c|\tau|.
\]
from which the proof follows.
\begin{flushright}
	\qed
\end{flushright}

%
In order to prove Theorem \ref{Thm:equivalence-GH-sums-systems} we need a result that follows from Theorem \ref{Thm:equivalence-gh}.
\begin{Prop}
	\label{Prop:charac-systems-coef-cte}
	The system $\mathcal X_0$ is globally hypoelliptic in $\TT^{n+1}$ if and only if the range of each $a_{10},\ldots,a_{n0}$ contains real numbers $\alpha_1,\ldots,\alpha_n$, respectively, that are not simultaneously approximable.
\end{Prop}
\begin{proof}
	Suppose the existence of such $\alpha_1,\ldots,\alpha_n$ and let $u \in \D'\left(\TT^{n+1}\right)$ such that $Xu \in \cinfty\left(\TT^{n+1}\right)$ for every $X \in \mathcal X_0$. In particular, if $L_j\doteq \partial_{t_j} + \alpha_j\partial_x$ for $j\in\{1,\ldots,n\}$ then $L_j u = f_j \in \cinfty\left(\TT^{n+1}\right)$ for every $j \in \{1,\ldots,n\}$ and Theorem \ref{Thm:equivalence-gh} shows that $u \in \cinfty\left(\TT^{n+1}\right)$.

	Conversely, suppose that the range of each $a_{10},\ldots,a_{n0}$ does not contain real numbers that are not simultaneously approximable. Since the non-Liouville numbers are dense in $\R$, it follows that each $a_{j0}$ is constant -- there is $\alpha_j \in \R$ such that $a_{j0}(t)=\alpha_j$ for every $t \in \TT^n$ -- and $\alpha_1,\ldots,\alpha_n$ is simultaneously approximable. Theorem \ref{Thm:equivalence-gh} (for $\mathcal X_0$) shows that $\mathcal X_0$ is not globally hypoelliptic in $\TT^{n+1}$.
\end{proof}
The proof of Theorem~\ref{Thm:equivalence-GH-sums-systems} follows from Theorem~\ref{Thm:equivalence-gh} and Proposition~\ref{Prop:charac-systems-coef-cte}.

\section{Global solvability and almost hypoellipticity}
\label{sec:global-solv}
We start by recalling some results of~\cite{afr22}: first, we equip $\cinfty\left(\TT^{N}\right)$ with its usual topology of Fréchet-Schwartz space induced by the Sobolev norms and say that an operator $Q:\cinfty\left(\TT^{N}\right)\to \cinfty\left(\TT^{N}\right)$ is globally solvable when $Q$ has closed range. Moreover, a system $\mathcal L$ of constant coefficients vector fields in $\TT^N$ is globally solvable when there exists a basis $\mathcal L_0=\left(L_1,\ldots,L_d\right)$ of $\Span_\R\mathcal L$ (as a subspace of the Lie agebra of $\TT^N$) such that the range of 
\[
\DD_{\mathcal L}:\cinfty\left(\TT^N\right)\to \cinfty\left(\TT^{N}\right)^d, \quad \DD_{\mathcal L}(u)=\left(L_1u,\ldots,L_du\right)
\]
is closed. This is equivalent to saying that the sublaplacian associated $\Delta_{\mathcal L} = -\sum_{j=1}^n L_j^2$ has closed range or that the following diophantine condition holds true: there exist $C,\rho>0$ such that 
\[
\max_{1\le j \le n} \left|L_j(\xi)\right| \ge C|\xi|^{-\rho}, \quad \xi \in \Z\setminus \Gamma,
\]
where
\[
 \Gamma\doteq \left\{\xi \in \Z^N : L_j(\xi) = 0,\,\, \forall j \in \{1,\ldots,n\}\right\}.
\]
When
\[
L_j = \partial_j + \alpha_j\partial_x, \quad \alpha_j \in \R, j \in \{1,\ldots,n\}
\]
in $\TT^{n+1}$, the diophantine condition is the following: there are $C,\rho>0$ such that
\[
\max_{1\le j \le m}\left|\tau_j+\alpha_j\xi\right| \ge C\left(1+|\xi|\right)^{-\rho}, \quad \forall \in (\tau,\xi) \in \Z^{n+1}\setminus \Gamma,
\]
where $\Gamma = \left\{(\tau,\xi) : \tau_j+\alpha_j\xi = 0, \,\, \forall j \in\{1,\ldots,n\}\right\}$. 
%
%
\begin{Thm}
	\label{Thm:charac-global-solv}
	The following properties are equivalent:
	\begin{enumerate}
		\item $P$ is globally solvable;
		\item the system $\mathcal X$ is globally solvable;
		\item the system $\mathcal X_0$ is globally solvable or there exists a point $t \in \TT^n$ of finite type for $X_1,\ldots,X_n$.
	\end{enumerate}
\end{Thm}
\begin{proof}
	We split our proof into cases. First we suppose that there is no point $t \in \TT^n$ of finite type for $X_1,\ldots,X_n$. In this case we take $S$ as in~\eqref{eq:simultaneous-reduction} and since $S$ conjugates $\mathcal X$ and $\mathcal X_0= \left\{Y_j\right\}_{j=1}^n$, it follows that $\mathcal X$ is globally solvable if and only if $\mathcal X_0$ is globally solvable. Furthermore, since
	\[
	\left(-\sum_{j=1}^n Y_j^2\right) S =  -\sum_{j=1}^n SX_j^2 = S P, 
	\]
	it follows that $P$ is globally solvable if and only if $-\sum_{j=1}^n Y_j^2$ is globally solvable, which in turn is equivalent (thanks to the considerations made before the theorem) to $\mathcal X_0$ being globally solvable.
	
	It remains to look at the case where there is a point of finite type for $X_1,\ldots,X_n$. Since every system of globally hypoelliptic operator (see~\cite{afjr23}) is also globally solvable, in this case the three items are true by Theorem~\ref{Thm:equivalence-gh} and there is nothing to prove.
\end{proof}
Now we recall that a system of operators $\mathcal S$ is said to be almost globally hypoelliptic when $u \in \D'\left(\TT^{n+1}\right)$ and $Xu \in \cinfty\left(\TT^{n+1}\right)$ for every $X\in S$ imply that there exists $v \in \cinfty\left(\TT^{n+1}\right)$ such that $X(u-v)=0$ for every $X \in \mathcal S$. It follows from~\cite{afjr23} that every finite almost globally hypoelliptic system is also globally solvable ($\mathcal S$ can be infinite if its operators have constant coefficients -- see~\cite{afr22}).
\\

\noindent{\bf Proof of Theorem~\ref{Thm:most-general-equiv}}: the strategy is similar to the proof of Theorem \ref{Thm:charac-global-solv}. Suppose that there is no point for finite type for $X_1,\ldots,X_n$. Then we use $S$ (recall~\eqref{eq:simultaneous-reduction}) to see that $P$ is globally solvable (respectively almost globally hypoelliptic) if and only if $P_0=-\sum_{j=1}^n Y_j^2$ is globally solvable (respectively almost globally hypoelliptic), which in turn is equivalent (see~\cite{afr22}) to the global solvability of $\mathcal X_0$ and also to its almost global hypoellipticity. By using $S$ again, $\mathcal X$ is globally solvable (respectively almost globally hypoelliptic) if and only if $\mathcal X_0$ is globally solvable (respectively almost globally hypoelliptic). So, in this case, all the items are equivalent.
	
	Now suppose that there is a point of finite type for $X_1,\ldots,X_n$. Then $P$ and $\mathcal X$ are globally hypoelliptic, which implies that $P$ and $\mathcal X$ are almost globally hypoelliptic and globally solvable, so all the items are true and there is nothing to prove.
\begin{flushright}
	\qed
\end{flushright}

\appendix

\section{Microlocal estimates}
\label{sec:me}

In this section we 	obtain microlocal estimates on the torus for elliptic operators by using the theory of pseudodifferential operators on torus presented in~\cite{rt_psdoas}. Given $u \in \D'\left(\TT^N\right)$, the set $\wf(u)$ denotes the frequencies $\xi \in \R^N\setminus\{0\}$ where $u$ is not smooth, that is $\xi_0 \notin \wf(u)$ if and only if there is an open cone $\Gamma\subset \R^N\setminus\{0\}$ containing $\xi_0$ such that for every $k\in\Z_+$ there exists $C_k>0$ such that
\[
\left|\hat u(\xi)\right| \le C_k\left(1+|\xi|\right)^{-k}, \quad \forall \xi \in \Gamma\cap \Z^N.
\]
%
%
%
Given $m\in\R$ we consider the class $S^m_{1,0}\left(\TT^N\times \Z^N\right)$ of the discrete symbols $a:\TT^N\times\Z^N\to \C$ such that $a(\cdot,\xi)$ is smooth for every $\xi \in \Z^N$ and satisfies the following inequality: given $\alpha,\beta \in \Z_+^N$, there exists $C>0$ such that
\[
\left|\Delta_\xi^\alpha\partial_x^\beta a(x,\xi)\right|\le C\left(1+|\xi|\right)^{m-|\alpha|},\quad \forall x \in \TT^N,\xi \in \Z^N,
\]
where $\Delta_\xi$ is a finite difference (see~\cite{rt_psdoas} for details). In particular, we have that given $k\in\Z_+$, there is $D_k>0$ such that
\begin{align}
	\label{Eq:decay-discrete-symbol}
	\left|\hat a(\eta,\xi-\eta)\right| \le D_k\left(1+|\eta|\right)^{-k}\left(1+|\xi-\eta|\right)^{m}, \quad \forall \xi,\eta \in \Z^N.
\end{align}
 The discrete symbol of every differential operator belongs to $S^m_{1,0}\left(\TT^N\times \Z^N\right)$. 

Every $a \in S^m_{1,0}\left(\TT^N\times\Z^N\right)$ defines a pseudodifferential operator 
\[
a(x,D) : \cinfty\left(\TT^N\right)\to \cinfty\left(\TT^N\right)
\]
by
\[
a(x,D)u = \sum_{\xi\in\Z^N} e^{ix\xi}a(x,\xi)\hat u(\xi), \quad \forall u \in \cinfty\left(\TT^N\right).
\]
In this case we say that $a(x,\xi)$ is the discrete symbol of $a(x,D)$. The reader can see~\cite{rt_psdoas} for the basic properties of these operators.

We also recall here Peetre's inequality:
\[
\left(1+|\zeta|\right)^\sigma \le \left(1+|\zeta'|\right)^{\sigma}\left(1+|\zeta-\zeta'|\right)^{|\sigma|}, \quad \forall \zeta,\zeta' \in \Z_+^N, \sigma \in \R.
\]
\begin{Prop}
	\label{Prop:inclusion_su}
	Let $a(x,D)$ be a linear partial pseudodifferential operator on $\TT^N$ of order $m$ and $u\in \D'\left(\mathbb T^N\right)$. Suppose that there exists an open cone $\Gamma\subset \R^N\setminus\{0\}$ such that for every $k\in\Z_+$ we can find $C_k>0$ satisfying
	\begin{equation}
		\label{inclusion_su_eq1}
		\left|\hat u(\xi)\right| \le C_k\left(1+|\xi|\right)^{-k}, \quad \forall \xi\in\Gamma\cap \Z^N.
	\end{equation}
	Then for any open cone $\Gamma'\subset \R^N\setminus\{0\}$ such that $\Gamma'\Subset \Gamma$, given $k\in \Z_+$ we can find a positive constant $C_k'>0$ such that
	\[
	\left|\widehat{\left(a(x,D)u\right)}(\xi)\right| \le C_k'\left(1+|\xi|\right)^{-k}, \quad \forall \xi\in\Gamma'\cap \Z^N.
	\]
	In particular, $\wf(a(x,D)u)\subset\wf(u)$.
\end{Prop}
\begin{proof}
	Choose $0<c<1$ such that $\eta \in \Gamma$ when $\xi\in\Gamma'$ and $|\xi-\eta|\le c|\xi|$. Then
	\begin{eqnarray*}
		\left|\widehat{\left(a(x,D)u\right)}(\xi)\right| &\le& \sum_{\eta \in \Z^N} \left|\hat a(\eta,\xi-\eta)\right| \left|\hat u(\xi-\eta)\right| \\
		&=& \sum_{c|\xi|\le |\eta|} \left|\hat a(\eta,\xi-\eta)\right| \left|\hat u(\xi-\eta)\right| + \sum_{|\eta|<c|\xi|} \left|\hat a(\eta,\xi-\eta)\right| \left|\hat u(\xi-\eta)\right|\\
		&\dot=& (I)+(II).
	\end{eqnarray*}
	\begin{enumerate}
		\item  Let us analyze the term $(I)$. Here we use only that $u$ belongs to some Sobolev space, that is there are $p\in\Z_+$ and $C>0$ such that
		$$
		|\hat u(\xi)|\le C(1+|\xi|)^p, \quad \forall \xi\in\Z^N.
		$$
		If $m' \in \Z_+$ satisfies $m'>m+p$ then for every $k\in\Z_+$ we use \eqref{Eq:decay-discrete-symbol} and Peetre's inequality to obtain
		\begin{eqnarray*}
			(I) &=& \sum_{c|\xi|\le |\eta|} \left|\hat a(\eta,\xi-\eta)\right| \left|\hat u(\xi-\eta)\right|\\
			&\le& \sum_{c|\xi|\le |\eta|} D_{k+2m'+N+1} (1+|\eta|)^{-k-2m'-N-1}(1+|\xi-\eta|)^m C (1+|\xi-\eta|)^p\\
			&\le& CD_{k+2m'+N+1}\sum_{c|\xi|\le |\eta|} (1+|\eta|)^{-k-2m'-N-1}(1+|\xi|)^{m+p}(1+|\eta|)^{m+p}\\
			&=& CD_{k+2m'+N+1}\sum_{c|\xi|\le |\eta|} (1+|\eta|)^{-k-2m'-N-1+m+p}(1+|\xi|)^{m+p}.
		\end{eqnarray*}
		Since we have $c|\xi|\le |\eta|$ in the last sum and $0<c<1$ it follows that
		$$
		1+|\xi| \le 1+c^{-1}|\eta| \le c^{-1}(1+|\eta|).
		$$
		The inequality $-m'+m+p\le 0$ yields
		\begin{eqnarray*}
			(I) &\le& CD_{k+2m'+N+1}\sum_{c|\xi|\le |\eta|} (1+|\eta|)^{-k-m'}(1+|\eta|)^{-N-1-m'+m+p}(1+|\xi|)^{m+p}\\
			&\le& CD_{k+m'+N+1}c^{-k-m'} \sum_{c|\xi|\le |\eta|} (1+|\xi|)^{-k-m'+m+p}(1+|\eta|)^{-N-1}.
		\end{eqnarray*}
		Hence, for a new constant $D_k'>0$, we conclude that
		\begin{equation}
			\label{inclusion_su_eq2}
			(I) \le \frac{D_k'}{(1+|\xi|)^k}, \quad \forall \xi\in\Gamma'.
		\end{equation}
		\item For the term $(II)$ we notice that by the construction of $c$ and $\Gamma'$ we can use \eqref{inclusion_su_eq1} with $\eta$ instead of $\xi$ if $\xi\in\Gamma'$ and $|\eta-\xi|\le c|\xi|$ since we have that $\eta \in \Gamma$ in this case. Take $m'\in\Z_+$ such that $m'>m$ and then for any $k,k'\in\Z_+$ we have
		\begin{eqnarray*}
			(II) &=& \sum_{|\eta|<c|\xi|} \left|\hat a(\eta,\xi-\eta)\right| \left|\hat u(\xi-\eta)\right|\\
			&\le& \left(\sup_{|\eta| < c|\xi|}\left|\hat u(\xi-\eta)\right|\right) \sum_{|\eta|<c|\xi|} \left|\hat a(\eta,\xi-\eta)\right|\\
			&=& \left(\sup_{|\xi-\eta| < c|\xi|}\left|\hat u(\eta)\right|\right) \sum_{|\eta|<c|\xi|} \left|\hat a(\eta,\xi-\eta)\right|\\
			&\le& \left(\sup_{|\xi-\eta| < c|\xi|} C_{k+m'}(1+|\eta|)^{-k-m'}\right) \sum_{|\eta|<c|\xi|} D_{k'}(1+|\eta|)^{-k'}(1+|\xi-\eta|)^m\\
			&\le& \left(\sup_{|\xi-\eta| < c|\xi|} C_{k+m'}(1+|\eta|)^{-k-m'}\right) \sum_{|\eta|<c|\xi|} D_{k'}(1+|\eta|)^{-k'+m}(1+|\xi|)^{m}.
		\end{eqnarray*}
		For $|\xi-\eta|<c|\xi|$ we have that
		$$
		|\xi| \le |\xi-\eta|+|\eta| < c|\xi|+|\eta|.
		$$
		So $(1-c)|\xi| < |\eta|$ and recalling that $0<c<1$ we obtain
		$$
		(1+|\eta|)^{-k-m'} \le (1+(1-c)|\xi|)^{-k-m'} \le (1-c)^{-k-m'}\left(1+|\xi|\right)^{-k-m'}.
		$$
		If we choose $k'$ such that $-k'+m < -N-1$ we have
		\begin{eqnarray*}
			(II) &\le& \sup_{|\xi-\eta| < c|\xi|} C_{k+m'}(1-c)^{-k-m'}(1+|\xi|)^{-k-m'+m} \sum_{|\eta|<c|\xi|} D_{k'}(1+|\eta|)^{-N-1}\\
			&\le& C_{k+m'}(1-c)^{-k-m'}(1+|\xi|)^{-k} \sum_{\eta\in\Z^N} D_{k'}(1+|\eta|)^{-N-1}, \quad \forall \xi\in\Gamma'.\\
		\end{eqnarray*}
		The proof follows from last inequality and \eqref{inclusion_su_eq2}.
	\end{enumerate}
\end{proof}
\begin{Lem}
	\label{symbol_vanishing_su}
	Let $a(x,\xi) \in S^m_{1,0}\left(\mathbb T^N\times\Z^N\right)$ such that $a(\cdot,\xi) = 0$ for every $\xi\in\Z^N\cap \Gamma$, where $\Gamma\subset \R^N\setminus\{0\}$ is an open cone containing a direction $\xi_0\in\R^N\setminus \{0\}$. If $u \in \D'\left(\mathbb T^N\right)$ then $\xi_0\notin \wf(a(x,D)u)$.
\end{Lem}
\begin{proof}
	Let $\Gamma'\Subset \Gamma$ be another open cone containing $\xi_0$ and $0<c<1$ such that $\eta \in \Gamma$ when $\xi\in\Gamma'$ and $|\eta-\xi|\le c |\xi|$. Then
	\begin{eqnarray*}
		\left|\widehat{a(x,D)u}(\xi)\right| &\le& \sum_{\eta \in \Z^N} \left|\hat a(\xi-\eta,\eta)\right| \left|\hat u(\eta)\right| \\
		&=& \sum_{|\xi-\eta|\ge c|\xi|} \left|\hat a(\xi-\eta,\eta)\right| \left|\hat u(\eta)\right|+\sum_{|\xi-\eta| < c|\xi|} \left|\hat a(\xi-\eta,\eta)\right| \left|\hat u(\eta)\right|\\
		&\dot=& (I)+(II), \quad \forall \xi\in\Gamma'\cap \Z^N.
	\end{eqnarray*}
	If $\xi\in\Gamma'$ and $|\eta-\xi| < c|\xi|$ then $a(\cdot,\eta)=0$ since $\eta \in \Gamma$. So we have
	$$
	\hat a(\xi-\eta,\eta) = (2\pi)^{-N}\int_{\mathbb T^N} e^{-i\left\langle x,\xi-\eta\right\rangle}a(x,\eta)dx = 0.
	$$
	Hence $(II)=0$ and we obtain that
	$$
	\left|\widehat{a(x,D)u}(\xi)\right| \le \sum_{|\xi-\eta| \ge c|\xi|} \left|\hat a(\xi-\eta,\eta)\right| \left|\hat u(\eta)\right| = \sum_{|\eta|\ge c|\xi|} \left|\hat a(\eta,\xi-\eta)\right| \left|\hat u(\xi-\eta)\right|.
	$$
	The right hand of this inequality is equal to the term (I) of Proposition \ref{Prop:inclusion_su}, which has a good decay independently of $u$. Thus we obtain \eqref{inclusion_su_eq2} in a similar way to that we have done in Proposition \ref{Prop:inclusion_su}.  
\end{proof}

For the proof of the next result the reader can see Theorem 4.5.3 of~\cite{rt_psdoas}.
\begin{Lem}
	\label{comparison_non_periodic}
	Suppose that $h_0=h_0(\xi)$ is a classical symbol of order $m$ and type $(1,0)$, that is $h_0 \in \cinfty\left(\mathbb R^N\right)$ and
	$$
	\left|D_\xi^\alpha h_0(\xi)\right| \le C_\alpha (1+|\xi|)^{m-|\alpha|}, \quad \forall \xi\in\R^N, \alpha\in\Z_+^N
	$$
	for some constants $C_\alpha>0$. If $h(x,\xi)=h(\xi)=h_0(\xi)$ for $(x,\xi)\in\TT^N\times\Z^N$, then $h \in S^m_{(1,0)}\left(\mathbb T^N\times\Z^N\right)$.
\end{Lem}

We obtain our microlocal estimates by looking at ``microlocal Sobolev spaces'', a notion precisely introduced in the sequence.
\begin{Def}
	\label{def_local_hs}
	Let $s\in\R$ and $\Gamma\subset\R^N\setminus \{0\}$ an open cone. We denote by $H^s_{\Gamma}\left(\mathbb T^N\right)$ the subspace of every $u\in \D'\left(\mathbb T^N\right)$ such that
	\[
	\left\|u\right\|_{s,\Gamma}^2 \doteq \sum_{\xi \in \Gamma\cap \Z^N} \left(1+|\xi|\right)^{2s}\left|\hat u(\xi)\right|^2 < \infty.
	\]
	Given $\xi \in \R^N\setminus\{0\}$, we denote by $H^s_{(\xi)}\left(\TT^N\right)$ the subspace of every $u\in\D'\left(\TT^N\right)$ such that there is an open cone $\Gamma\subset\R^N\setminus\{0\}$ such that $\xi\in\Gamma$ and $u\in H^s_{\Gamma}\left(\TT^N\right)$.
\end{Def}
\begin{Prop}
	\label{relation_local_hs}
	Let $u\in \D'\left(\mathbb T^N\right)$ and $\Gamma\subset \R^N\setminus\{0\}$ an open cone. Then $u \in H_{(\xi)}^s\left(\mathbb T^N\right)$ for every $\xi\in\Gamma$ if and only if $u\in H^s_{\Gamma'}\left(\TT^N\right)$ for any open cone $\Gamma'\Subset \Gamma$.
\end{Prop}
\begin{proof}
Suppose that $u\in 	H_{(\xi)}^s\left(\mathbb T^N\right)$ for every $\xi\in\Gamma$ and let $\Gamma'\Subset\Gamma$ be an open cone. Given $\xi\in\overline{\Gamma'}\cap S^{N-1}$, since $\xi\in\Gamma$, there is an open cone $\Gamma_\xi\subset\R^N\setminus\{0\}$ such that $u \in H^s_{\Gamma_\xi}\left(\TT^N\right)$. By the compactness of $\overline{\Gamma'}\cap S^{N-1}$, there are $\xi_1,\ldots,\xi_m\in \Gamma'\cap S^{N-1}$ such that $\Gamma'\subset \Gamma_{\xi_1}\cup\ldots\cup\Gamma_{\xi_m}$. Thus
\begin{align*}
	\left\|u\right\|_{s,\Gamma'}^2 &= \sum_{\xi\in\Gamma'\cap \Z^N} \left(1+|\xi|\right)^{2s}\left|\hat u(\xi)\right|^2 \le \sum_{j=1}^m \left\|u\right\|_{s,\Gamma_{\xi_j}}^2<\infty.
\end{align*}
The converse is trivial since every $\xi \in \Gamma$ belongs to some $\Gamma'\Subset\Gamma$.
\end{proof}

\begin{Prop}
	\label{local_hs_su}
	Let $u\in \D'\left(\mathbb T^N\right)$ and $\xi_0\in\R^N\setminus \{0\}$. Then $\xi_0 \notin \wf(u)$ if and only if there exist an open cone $\Gamma\subset\R^N\setminus\{0\}$ such that $\xi_0\in\Gamma$ and $u \in H_{\Gamma}^s\left(\mathbb T^N\right)$ for every $s\in\R$.
\end{Prop}
\begin{proof}
	If $\xi_0 \notin \wf(u)$ then we can find an open cone $\Gamma\subset \R^N\setminus\{0\}$, with $\xi_0 \in \Gamma$, such that given $k\in\Z_+$ there exists $C_k>0$ such that
	\[
		\left|\hat u(\eta)\right| \le \frac{C_k}{(1+|\eta|)^k}, \quad \forall \eta \in \Gamma\cap \Z^N.
	\]
	Given $s\in\R$, we choose $k\in \N$ such that $s+N<k$ to conclude that
	\begin{align*}
		\left\|u\right\|_{s,\Gamma}^2 &= \sum_{\xi\in\Gamma\cap\Z^N} \left(1+|\xi|\right)^{2s}\left|\hat u(\xi)\right|^2\\
		&\le C_{k}\sum_{\xi\in\Gamma} \left(1+|\xi|\right)^{2s}\left(1+|\xi|\right)^{-2k}\\
		&\le C_{k}\sum_{\xi\in\Gamma} \left(1+|\xi|\right)^{-2N}\\
		&<\infty.
	\end{align*}
	
	For the converse we suppose that there is an open cone $\Gamma\subset\R^{N}\setminus\{0\}$ containing $\xi_0$ and such that $u \in H_{\Gamma}^s\left(\mathbb T^N\right)$ for every $s\in\R$. Then given $k\in\Z_+$, we obtain that
	\begin{align*}
		\left(1+|\xi|\right)^{2k}\left|\hat u(\xi)\right|^2 \le \left\|u\right\|^2_{k,\Gamma},\quad \forall \xi\in\Gamma\cap \Z^N,
	\end{align*}
	from which follows that $\xi_0\notin \mathfrak s(u)$.
\end{proof}

Now we consider a definition of ellipticity for our class of pseudodifferential operators in a fixed direction (see Definition 4.9.1 of~\cite{rt_psdoas}):
\begin{Def}
	\label{local_ellipticity}
	Let $a(x,\xi) \in S^m_{1,0}\left(\mathbb T^N\times\Z^N\right)$ and $\xi_0 \in \R^N\setminus \{0\}$. We say that $a(x,D)$ (or $a(x,\xi)$) is elliptic (of order $m$) in $\xi_0$ if there exist an open cone $\Gamma$ containing $\xi_0$ and constants $C,r>0$ such that
	$$
	|a(x,\xi)| \ge C(1+|\xi|)^m, \quad \forall x\in\TT^N, \xi\in\Gamma\cap \Z^N, |\xi|\ge r.
	$$
\end{Def}
Let $a \in S^m_{1,0}\left(\mathbb T^N\times \Z^N\right)$ and consider $a(x,D)$ the operator associated to $a$. A principal symbol of $a(x,D)$ is a symbol $b \in S^m_{1,0}\left(\mathbb T^N\times \Z^N\right)$ such that $a-b \in S^{m-1}_{1,0}\left(\mathbb T^N\times \Z^N\right)$. The next result shows that the ellipticity is a property about principal symbols.
\begin{Prop}
	\label{local_ellipticity_principal_symbol}
	Let $a(x,\xi) \in S^m_{1,0}\left(\mathbb T^N\right)$ and $\xi_0 \in \R^N\setminus \{0\}$. Then $a(x,D)$ is elliptic in $\xi_0$ if and only if $b(x,\xi)$ is elliptic in $\xi_0$, where $b(x,\xi)$ is any principal symbol of $a(x,D)$. 
\end{Prop}
\begin{proof}
	
	Suppose that $a(x,\xi)$ is elliptic in $\xi_0$. Then there exist an open cone $\Gamma$ containing $\xi_0$ and constants $C,r>0$ such that
	$$
	|a(x,\xi)| \ge C(1+|\xi|)^m, \quad \forall \xi\in\Gamma\cap \Z^N, |\xi|\ge r.
	$$
	If $c(x,\xi) = a(x,\xi)-b(x,\xi)$, where $b(x,\xi)$ is a principal symbol of $a(x,D)$, then $c(x,\xi) \in S^{m-1}\left(\mathbb T^N\right)$ and we can find $D>0$ such that
	$$
	\left|c(x,\xi)\right| \le D(1+|\xi|)^{m-1}, \quad \forall \xi\in\Z^N.
	$$
	Then, for $\xi\in\Gamma, |\xi| \ge r$, we obtain
	\begin{eqnarray*}
		\left|b(x,\xi)\right| &=& \left|a(x,\xi)-c(x,\xi)\right|\\
		&\ge& \left|a(x,\xi)\right| - \left|c(x,\xi)\right|\\
		&\ge& C(1+|\xi|)^m - D(1+|\xi|)^{m-1}\\
		&=& (1+|\xi|)^m\left(C-\frac{D}{(1+|\xi|)}\right).
	\end{eqnarray*}
	So if $\rho>0$ is such that
	$$
	(1+|\xi|) \ge \frac{2D}{C}, \ \mbox{if} \ |\xi|>\rho
	$$
	then we conclude that
	$$
	\left|b(x,\xi)\right| \ge \frac{C}{2}(1+|\xi|)^m
	$$
	for all $\xi\in\Gamma$ such that $|\xi| > \max\{r,\rho\}$. 
	
	For the converse we suppose that $b(x,\xi)$ is a principal symbol of $a(x,D)$ that is elliptic in $\xi_0$. We consider again $c(x,\xi) = a(x,\xi)-b(x,\xi) \in S^{m-1}\left(\mathbb T^N\times \Z^N\right)$ and we can repeat the argument used above with the roles of $a(x,\xi)$ and $b(x,\xi)$ changed to conclude that $a(x,\xi)$ is elliptic in $\xi_0$. 
\end{proof}
%
%
%
\begin{Rem}
	Let $s \in \R, \xi_0\in\R^N\setminus\{0\}$ and $u\in \D'\left(\TT^N\right)$. In the next result we shall use the following characterization: $u\in H^s_{(\xi_0)}\left(\TT^N\right)$ if and only if $u=v+w$, where $v \in H^s\left(\TT^N\right)$ and $w \in \D'\left(\TT^N\right)$ is such that $\xi_0\notin \mathfrak s(w)$, where $H^s\left(\TT^N\right)$ stands for the usual Sobolev spaces on $\TT^N$. Indeed, if $u \in H^s_{(\xi_0)}\left(\TT^N\right)$ then there is an open cone $\Gamma\subset\R^N\setminus\{0\}$ containing $\xi_0$ and such that $\left\|u\right\|_{s,\Gamma}<\infty$. If we set
	\[
	\hat v(\xi)
	=
	\begin{cases}
		\hat u(\xi),&\quad \xi \in \Gamma\cap \Z^N\\
		0,&\quad \mbox{otherwise}
	\end{cases}
	\quad\mbox{and}\quad
		\hat w(\xi)
	=
	\begin{cases}
		0,&\quad \xi \in \Gamma\cap \Z^N\\
		\hat u(\xi),&\quad \mbox{otherwise}
	\end{cases}
	\]
	then $u=v+w$, $v\in H^s\left(\TT^N\right)$ and $w \in \D'\left(\TT^N\right)$ satisfy $\xi_0\notin \mathfrak s(w)$. For the converse it suffices to note that $H^s\left(\TT^N\right)\subset H^s_\Gamma\left(\TT^N\right)$ for any open cone $\Gamma\subset\R^N\setminus\{0\}$ and use Proposition~\ref{local_hs_su}.
\end{Rem}
\begin{Thm}
	\label{elipticity_local_hs}
	Let $a(x,\xi) \in S^m_{1,0}\left(\mathbb T^N\times \Z^N\right)$ such that $a(x,D)$ is elliptic in $\xi_0 \in \R^N\setminus \{0\}$. If $u\in\D'\left(\TT^N\right)$ and $a(x,D)u \in H_{(\xi_0)}^s\left(\mathbb T^N\right)$ then $u \in H^{s+m}_{(\xi_0)}\left(\mathbb T^N\right)$.
\end{Thm}
\begin{proof}
	
	We can assume, without loss of generality, that $a(x,\xi)\ge 0$ for all $(x,\xi)\in\mathbb T^N\times\Z^N$. Indeed, let $b(x,D) = a^\ast(x,D) a(x,D) \in S^{2m}_{1,0}\left(\mathbb T^N\times\Z^N\right)$. By Proposition \ref{local_ellipticity_principal_symbol}, since $|a(x,\xi)|^2$ is a principal symbol of $b(x,D)$, we have that $b(x,D)$ is elliptic (of order $2m$) in $\xi_0$. Furthermore $b(x,D)u \in H^{s-m}_{(\xi_0)}\left(\mathbb T^N\right)$: $a(x,D)u=v+w$ with $v \in H^s\left(\mathbb T^N\right)$ and $w \in \D'\left(\mathbb T^N\right)$ satisfies $\xi_0 \notin \wf(w)$. Hence
	$$
	a^\ast(x,D)a(x,D)u = a^\ast(x,D)v + a^\ast (x,D) w, 
	$$ 
	with $a^\ast (x,D) v \in H^{s-m}\left(\mathbb T^N\right)$ and $\xi_0 \notin \wf(a^\ast(x,D) w)$ by Proposition \ref{Prop:inclusion_su}. So from now on we assume that $a(x,\xi) \ge 0$ for all $x \in \mathbb T^N$ and $\xi\in\Z^N$.
	
	It follows from the ellipticity of $a(x,\xi)$ in $\xi_0$ that there exist an open cone $\Gamma\subset\R^N\setminus\{0\}$ containing $\xi_0$ and constants $C,r>0$ such that
	$$
	a(x,\xi) \ge C(1+|\xi|)^m, \quad \forall \xi\in\Gamma\cap \Z^N, |\xi|\ge r.
	$$
	Let $\Gamma_0\subset\R^N\setminus\{0\}$ be an open cone containing $\xi_0$ such that $\Gamma_0 \Subset \Gamma$ e $h \in \cinfty\left(\mathbb \R^N\right)$ such that $\supp h \subset \Gamma, h\equiv 1$ in $\Gamma_0$ and $h$ satisfies the condition of Lemma \ref{comparison_non_periodic} for $m=0$. We define
	$$
	a'(x,\xi) = h(\xi)a(x,\xi) + (1-h(\xi))(1+|\xi|^2)^{m/2}, \quad \forall x \in \mathbb T^N, \xi\in\Z^N.
	$$
	Then $a'(x,\xi) \in S_{1,0}^m\left(\mathbb T^N\times\Z^N\right)$ is elliptic. Indeed, we first increase $r>0$ in such way that $(1+|\xi|^2)^{m/2} \ge D(1+|\xi|)^m$ for all $|\xi|\ge r$ and some $D>0$. If $\xi\in \Z^N\setminus\Gamma$ then $h(\xi)=0$ and 
	$$
	a'(x,\xi) = (1+|\xi|^2)^{m/2} \ge D(1+|\xi|)^m.
	$$ 
	If $\xi\in\Gamma$ and $|\xi|\ge r$ then
	\begin{eqnarray*}
		a'(x,\xi) &\ge& C h(\xi)(1+|\xi|)^m + D(1-h(\xi))(1+|\xi|)^m\\
		&\ge& \min\{D,C\}\left(h(\xi)(1+|\xi|)^m + (1-h(\xi))(1+|\xi|)^m\right)\\
		&=& \min\{D,C\}(1+|\xi|)^m. 
	\end{eqnarray*}
	If we denote $A'=a'(x,D)$, it follows from Theorem 4.9.6 of~\cite{rt_psdoas} that there exists a parametrix $B$ for $A'$, that is an operator $B=b(x,D)$ for some $b \in S^{-m}_{1,0}\left(\TT^N\times\Z^N\right)$ such that the operators $A'B-I,BA'-I$ map $\D'\left(\TT^N\right)$ into $\cinfty\left(\TT^N\right)$. If $A=a(x,D)$ then by hypothesis we can write $Au = v + w$, where $v \in H^s\left(\TT^N\right)$ and $\xi_0\notin \wf(w)$, from which we obtain that
	\begin{align*}
		u &= \left(I-BA'\right)u + BA'u \\
		&= \left(I-BA'\right)u + B(A'-A)u + BAu\\
		&= Bv + \left(I-BA'\right)u + B(A'-A)u + Bw.
	\end{align*}
	Note that  $\xi_0 \notin \wf\left(\left(I-BA'\right)u + B(A'-A)u + Bw\right)$: $(I-BA')u \in \cinfty\left(\TT^N\right)$, $\xi_0\notin \mathfrak s(Bw)$ by Proposition \ref{Prop:inclusion_su} and since $a(x,\xi)-a'(x,\xi)=0$ for all $x\in\mathbb T^N$ and $\xi\in\Gamma_0$, it follows from Lemma \ref{symbol_vanishing_su} and from Proposition \ref{Prop:inclusion_su} that $\xi_0 \notin \wf\left(B(A'-A)u\right)$. The proof follows since $Bv \in H^{s+m}\left(\TT^N\right)$.
\end{proof}

\begin{Cor}
	\label{local_ellipticity_su}
	Let $a(x,\xi) \in S^m_{1,0}\left(\mathbb T^N\times\Z^N\right)$ be an elliptic symbol in $\xi_0 \in \R^N\setminus \{0\}$. If $u \in \D'\left(\mathbb T^N\right)$ and $\xi_0 \notin \wf(a(x,D)u)$ then $\xi_0 \notin \wf(u)$.
\end{Cor}
\begin{proof}
	Since $\xi_0 \notin \wf(a(x,D)u)$ we can use Proposition \ref{local_hs_su} to obtain an open cone $\Gamma\subset\R^N\setminus \{0\}$ containing $\xi_0$ such that $a(x,D) u \in H^s_{\Gamma}\left(\mathbb T^N\right)$ for every $s\in\R$. Since $a(x,D)$ is elliptic in $\xi_0$, we can shrink $\Gamma$ if necessary to suppose that $a(x,\xi)$ is elliptic in every $\xi \in \Gamma$. Since $a(x,D)u \in H^s_{\Gamma}\left(\TT^N\right)$, it follows by definition that $a(x,D)u \in H^s_{(\xi)}\left(\TT^N\right)$ for every $\xi \in \Gamma$. By Theorem~\ref{elipticity_local_hs} we obtain that $u \in H^{s+m}_{(\xi)}\left(\TT^N\right)$ for every $\xi \in \Gamma$ and $s \in \R$. So if we fix $\Gamma'\Subset \Gamma$, we have that $u \in H^{s+m}_{\Gamma'}\left(\TT^N\right)$ for every $s \in \R$ by Proposition~\ref{relation_local_hs}. Since $s$ is arbitrary, we use again Proposition~\ref{local_hs_su} to conclude that $\xi_0\notin \mathfrak s(u)$.
\end{proof}

We finish this section applying Corollary \ref{local_ellipticity_su} to a specific class of operators that depend only on the $t$-variable in a product $\TT^n_t\times \TT^m_x$.
%
The proof of the next result is similar to Theorem 3.3 of~\cite{cc17} or Lemma 4.2 of~\cite{fp24}.
\begin{Cor}
	\label{Cor:ellipticity-wf}
	Let $P=P\left(t,x,D_t,D_x\right)$ be a linear partial differential operator on $\TT_t^n\times\TT_x^m$ and suppose that $P$ is elliptic in any point $(t,x,\tau,0)\in \TT^{n+m}\times\left(\Z^{n+m}\setminus\{0\}\right)$. If $u \in \D'\left(\TT^n\times \TT^m\right)$ and $Pu \in \cinfty\left(\TT^n\times\TT^m\right)$ then there exists $c>0$ such that for every $k\in\Z_+$ there is $C_k>0$ such that
	\[
	\left|\hat u(\tau,\xi)\right| \le C_k\left(1+|(\tau,\xi)|\right)^{-k}, \quad |\xi| \le c|\tau|.
	\]
\end{Cor}

\bibliographystyle{plain}
\bibliography{bibliography}
\end{document}